\documentclass{amsart}
\usepackage{amsopn}
\usepackage{amssymb, amscd}

\newcommand{\nc}{\newcommand}

\nc{\vg}{\mathfrak{v} } \nc{\wg}{\mathfrak{w} } \nc{\zg}{\mathfrak{z} }
\nc{\ngo}{\mathfrak{n} } \nc{\kg}{\mathfrak{k} } \nc{\mg}{\mathfrak{m} }
\nc{\bg}{\mathfrak{b} } \nc{\ggo}{\mathfrak{g} } \nc{\ggob}{\overline{\mathfrak{g}}
} \nc{\sog}{\mathfrak{so} } \nc{\sug}{\mathfrak{su} } \nc{\spg}{\mathfrak{sp} }
\nc{\slg}{\mathfrak{sl} } \nc{\glg}{\mathfrak{gl} } \nc{\cg}{\mathfrak{c} }
\nc{\rg}{\mathfrak{r} } \nc{\hg}{\mathfrak{h} } \nc{\tg}{\mathfrak{t} }
\nc{\ug}{\mathfrak{u} } \nc{\dg}{\mathfrak{d} } \nc{\ag}{\mathfrak{a} }
\nc{\pg}{\mathfrak{p} } \nc{\sg}{\mathfrak{s} } \nc{\pca}{\mathcal{P}}
\nc{\nca}{\mathcal{N}} \nc{\lca}{\mathcal{L}} \nc{\oca}{\mathcal{O}}
\nc{\mca}{\mathcal{M}} \nc{\tca}{\mathcal{T}} \nc{\aca}{\mathcal{A}}
\nc{\cca}{\mathcal{C}} \nc{\gca}{\mathcal{G}} \nc{\sca}{\mathcal{S}}
\nc{\hca}{\mathcal{H}} \nc{\bca}{\mathcal{B}} \nc{\dca}{\mathcal{D}}
\nc{\val}{\operatorname{val}}

\nc{\vp}{\varphi} \nc{\ddt}{\tfrac{{\rm d}}{{\rm d}t}} \nc{\im}{\mathtt{i}}

\nc{\SO}{\mathrm{SO}} \nc{\Spe}{\mathrm{Sp}} \nc{\Sl}{\mathrm{SL}}
\nc{\SU}{\mathrm{SU}} \nc{\Or}{\mathrm{O}} \nc{\U}{\mathrm{U}} \nc{\Gl}{\mathrm{GL}}
\nc{\Se}{\mathrm{S}} \nc{\Cl}{\mathrm{Cl}} \nc{\Spein}{\mathrm{Spin}}
\nc{\Pin}{\mathrm{Pin}} \nc{\G}{\mathrm{GL}_n} \nc{\g}{\mathfrak{gl}_n}

\nc{\RR}{{\Bbb R}} \nc{\HH}{{\Bbb H}} \nc{\CC}{{\Bbb C}} \nc{\ZZ}{{\Bbb Z}}
\nc{\FF}{{\Bbb F}} \nc{\NN}{{\Bbb N}} \nc{\QQ}{{\Bbb Q}} \nc{\PP}{{\Bbb P}}

\nc{\vs}{\vspace{.2cm}} \nc{\vsp}{\vspace{1cm}} \nc{\ip}{\langle\cdot,\cdot\rangle}
\nc{\ipp}{(\cdot,\cdot)} \nc{\la}{\langle} \nc{\ra}{\rangle} \nc{\unm}{\tfrac{1}{2}}
\nc{\unc}{\tfrac{1}{4}} \nc{\und}{\tfrac{1}{16}} \nc{\no}{\vs\noindent}
\nc{\lam}{\Lambda^2(\RR^n)^*\otimes\RR^n} \nc{\tangz}{{\rm T}^{\rm Zar}}
\nc{\nor}{{\sf n}} \nc{\eigen}{(k_1<...<k_r;d_1,...,d_r)}
\nc{\eigencero}{(0<k_2<...<k_r;d_1,...,d_r)} \nc{\mum}{/\!\!/}
\nc{\kir}{/\!\!/\!\!/} \nc{\Ri}{\tfrac{4}{||\mu||^2}\Ric_{\mu}}
\nc{\ds}{\displaystyle} \nc{\ben}{\begin{enumerate}} \nc{\een}{\end{enumerate}}
\nc{\f}{\frac}

\nc{\He}{\operatorname{Hess}} \nc{\ad}{\operatorname{ad}}
\nc{\Ad}{\operatorname{Ad}} \nc{\rank}{\operatorname{rank}}
\nc{\Irr}{\operatorname{Irr}} \nc{\End}{\operatorname{End}}
\nc{\Aut}{\operatorname{Aut}} \nc{\Inn}{\operatorname{Inn}}
\nc{\Der}{\operatorname{Der}} \nc{\Ker}{\operatorname{Ker}}
\nc{\Iso}{\operatorname{I}} \nc{\Diff}{\operatorname{D}} \nc{\Lie}{\operatorname{L}}
\nc{\tr}{\operatorname{tr}} \nc{\dif}{\operatorname{d}}
\nc{\sen}{\operatorname{sen}} \nc{\modu}{\operatorname{mod}}
\nc{\Ric}{\operatorname{R}} \nc{\Ricg}{\operatorname{Ric^{\gamma}}}
\nc{\Ricc}{\operatorname{Ric^{c}}} \nc{\sym}{\operatorname{sym}}
\nc{\symac}{\operatorname{sym^{ac}}} \nc{\symc}{\operatorname{sym^{c}}}
\nc{\scalar}{\operatorname{sc}} \nc{\grad}{\operatorname{grad}}
\nc{\ricci}{\operatorname{ric}} \nc{\ricciac}{\operatorname{ric^{ac}}}
\nc{\riccic}{\operatorname{ric^{c}}} \nc{\riccig}{\operatorname{ric^{\gamma}}}
\nc{\Rin}{\operatorname{M}} \nc{\Le}{\operatorname{L}} \nc{\tang}{\operatorname{T}}
\nc{\level}{\operatorname{level}} \nc{\rad}{\operatorname{r}}
\nc{\abel}{\operatorname{ab}} \nc{\CH}{\operatorname{CH}}
\nc{\mcc}{\operatorname{mcc}} \nc{\Adj}{\operatorname{Adj}}

\theoremstyle{plain}
\newtheorem{theorem}{Theorem}[section]
\newtheorem{proposition}[theorem]{Proposition}
\newtheorem{corollary}[theorem]{Corollary}
\newtheorem{lemma}[theorem]{Lemma}

\theoremstyle{definition}
\newtheorem{definition}[theorem]{Definition}

\theoremstyle{remark}
\newtheorem{remark}[theorem]{Remark}

\newtheorem{example}[theorem]{Example}

\title{Einstein solvmanifolds: existence and non-existence questions}

\author{Jorge Lauret} \author{Cynthia Will}

\address{FaMAF and CIEM, Universidad Nacional de C\'ordoba, 5000 C\'ordoba, Argentina}
\email{lauret@famaf.unc.edu.ar, cwill@famaf.unc.edu.ar}

\thanks{
Supported by grants from Fundaci\'on Antorchas, CONICET, FONCyT and SeCyT (Univ. Nac. de C\'ordoba)}

\begin{document}

\maketitle

\begin{abstract}
The aim of this paper is to study the problem of which solvable Lie groups admit an
Einstein left invariant metric. The space $\nca$ of all nilpotent Lie brackets on
$\RR^n$ parametrizes a set of $(n+1)$-dimensional rank-one solvmanifolds $\{
S_{\mu}:\mu\in\nca\}$, containing the set of all those which are Einstein in that
dimension. The moment map for the natural $\G$-action on $\nca$, evaluated at
$\mu\in\nca$, encodes geometric information on $S_{\mu}$ and suggests the use of
strong results from geometric invariant theory. For instance, the functional on
$\nca$ whose critical points are precisely the Einstein $S_{\mu}$'s, is the square
norm of this moment map. We use a $\G$-invariant stratification for the space
$\nca$ and show that there is a strong interplay between the strata and the Einstein
condition on the solvmanifolds $S_{\mu}$.  As an application, we obtain criteria to
decide whether a given nilpotent Lie algebra can be the nilradical of a rank-one
Einstein solvmanifold or not.  We find several examples of $\NN$-graded (even
$2$-step) nilpotent Lie algebras which are not.  A classification in the
$7$-dimensional, $6$-step case and an existence result for certain $2$-step algebras
associated to graphs are also given.
\end{abstract}

\section{Introduction}\label{intro}

The study of homogeneous Einstein manifolds breaks into two very distinguishable
parts, compact and noncompact cases, according to the sign of the Einstein constant.
Known and expected results as well as approaches and techniques are quite different
in one and other case, although they both share the same basic general question:

\begin{quote}
{\bf Problem 1}.  Which homogeneous spaces $G/K$ admit a $G$-invariant Einstein
metric?
\end{quote}

In the noncompact case, the only known examples until now are all of a very
particular kind; namely, simply connected solvable Lie groups endowed with a left
invariant metric $(S,\ip)$ ({\it solvmanifolds}). According to the long standing
Alekseevskii conjecture (see \cite[7.57]{Bss}), every noncompact homogeneous Einstein manifold should be of this kind.  It has recently been
proved in \cite{standard} that any Einstein solvmanifold is {\it standard}: if
$\sg=\ag\oplus\ngo$ is the orthogonal decomposition of the Lie algebra $\sg$ of $S$
with $\ngo=[\sg,\sg]$, then $[\ag,\ag]=0$.  Standard Einstein solvmanifolds were
extensively studied by J. Heber in \cite{Hbr}, who proved very nice structural and
uniquenes results. Since a solvable Lie group admits at most one standard Einstein
metric up to isometry and scaling (see \cite[Theorem 5.1]{Hbr}), a substantial part
of Problem 1 in the noncompact case (probably the whole of it) can be reformulated
as

\begin{quote}
{\bf Problem 2}.  Which solvable Lie groups admit a standard Einstein metric?
\end{quote}

To approach this problem is the goal of the present paper.  We first recall that the
study of standard Einstein solvmanifolds reduces to the rank-one case, that is,
$\dim{\ag}=1$ (see \cite[Sections 4.5,4.6]{Hbr}). Our second observation is that
actually everything is determined by the nilpotent Lie algebra $\ngo=[\sg,\sg]$.
Indeed, a nilpotent Lie algebra $\ngo$ is the nilradical of a rank-one Einstein
solvmanifold if and only if $\ngo$ admits a {\it nilsoliton} metric (also called a {\it minimal} metric), meaning that its Ricci operator is a multiple of the identity
modulo a derivation of $\ngo$.  We call such an $\ngo$ an {\it Einstein nilradical}.
A nilsoliton metric on a nilpotent Lie algebra is also known to be unique up to
isometry and scaling (see \cite[Theorem 3.5]{soliton}), and so we can reformulate
Problem 2 as the equivalent

\begin{quote}
{\bf Problem 3}.  Which nilpotent Lie algebras are Einstein nilradicals?
\end{quote}

Any nilpotent Lie algebra of dimension at most $6$ is an Einstein nilradical (see
\cite{finding, Wll}), and the only known obstruction up to now is that any Einstein
nilradical $\ngo$ has to admit an $\NN$-gradation, that is,
$\ngo=\ngo_1\oplus...\oplus\ngo_r$, $[\ngo_i,\ngo_j]\subset\ngo_{i+j}$ (see
\cite[Theorem 4.14]{Hbr}).  The $\NN$-gradation comes from the eigenspace
decomposition for the derivation defined by $\ag$ as its eigenvalues form a set of
natural numbers without a common divisor called the {\it eigenvalue type} of the
standard Einstein solvmanifold.  A new natural question then arises:

\begin{quote}
{\bf Problem 4}.  Is every $\NN$-graded nilpotent Lie algebra an Einstein
nilradical?
\end{quote}

We shall give a negative answer to this problem by exhibiting several
counterexamples with many different features.  This will show in particular the high
difficulty of Problem 3.

The following approach allows us to use several results from geometric invariant
theory.  Fix an inner product vector space
$$
(\sg=\RR H\oplus\RR^n,\ip),\qquad \la H,\RR^n\ra=0,\quad \la H,H\ra=1.
$$
For each nilpotent Lie bracket $\mu:\RR^n\times\RR^n\longrightarrow\RR^n$, there
exists a unique derivation $D_{\mu}\in\Der(\mu)$ such that the $(n+1)$-dimensional
rank-one solvmanifold $S_{\mu}$ modelled on $(\sg=\RR H\oplus\RR^n,\ip)$ by
$$
[H,X]=D_{\mu}X,\qquad [X,Y]=\mu(X,Y), \qquad X,Y\in\RR^n,
$$
has a chance of being Einstein (see \cite{critical} or Section \ref{pre}, where a
background of Einstein solvmanifolds is given). Thus the space $\nca$ of all
nilpotent Lie brackets on $\RR^n$ parametrizes a set of $(n+1)$-dimensional rank-one
solvmanifolds
$$
\{ S_{\mu}:\mu\in\nca\},
$$
that contains the much smaller set of all those which are Einstein in that
dimension. $\nca$ is an algebraic subset of the vector space
$V=\Lambda^2(\RR^n)^*\otimes\RR^n$, which is invariant with respect to the change of
basis action of $\G$ on $V$ and Lie algebra isomorphism classes correspond to
$\G$-orbits. Other than the solvmanifold $S_{\mu}$, there is another Riemannian
manifold naturally associated to each $\mu\in\nca$; namely, the nilmanifold
$(N_{\mu},\ip)$, where $N_{\mu}$ is the simply connected nilpotent Lie group with
Lie algebra $(\RR^n,\mu)$ endowed with the left invariant metric determined by
$\ip|_{\RR^n\times\RR^n}$. The orbit $\G.\mu$ can be viewed as a parametrization of
the set of all left invariant metrics on $N_{\mu}$, and hence

\begin{quote}
$\mu$ is an Einstein nilradical if and only if $S_{g.\mu}$ is Einstein for some
$g\in\G$.
\end{quote}

The set $\G.\mu$ is in general huge, even up to the $\Or(n)$-action (or
equivalently up to isometry when viewed as a set of Riemannian metrics), and here is
where one finds the main difficulty in trying to decide whether a given nilpotent
Lie algebra is an Einstein nilradical or not. Such a problem is considerably easier
to tackle if one knows a priori the eigenvalue type.

In geometric invariant theory, a moment map for linear actions of complex reductive
Lie groups has been defined in \cite{Nss} and \cite{Krw1}.  A remarkable fact is
that, in our situation (i.e. for the $\G$-action on $V$), this moment map
$m:V\longrightarrow\glg_n$ encodes geometric information on $S_{\mu}$; indeed, it
was proved in \cite{minimal} that
$$
m(\mu)=\tfrac{4}{||\mu||^2}\Ric_{\mu}, \qquad\forall\mu\in\nca,
$$
where $\Ric_{\mu}$ is the Ricci operator of $(N_{\mu},\ip)$.  We also know that
$S_{\mu}$ is Einstein if and only if $\mu$ is a critical point of the functional
$F:V\longrightarrow\RR$ defined by
$$
F(\mu)=\frac{16\tr{\Ric_{\mu}^2}}{||\mu||^4}
$$
(see \cite{critical}).  Thus $F$ is precisely the square norm of the moment map, a
function extensively studied in symplectic geometry and with very nice properties of
convexity and minimality relative to the orbit structure of the action.  Some of the
structural and uniqueness results on standard Einstein solvmanifolds obtained in
\cite{Hbr} follow indeed as applications of such properties (see \cite{praga}).

In Section \ref{st}, we consider a $\G$-invariant stratification for the vector
space $V$ with certain frontier properties defined in \cite{standard}.  The strata
are parametrized by a finite set $\bca$ of diagonal $n\times n$ matrices.  We prove
that if $\mu\in\nca$ belongs to the stratum $\sca_{\beta}$, $\beta\in\bca$, then
$F(\mu)\geq ||\beta||^2$ and equality holds if and only if $S_{\mu}$ is Einstein, if
and only if $\tfrac{4}{||\mu||^2}\Ric_{\mu}$ is conjugate to $\beta$. Conversely, if
$S_{\mu}$ is Einstein then $\mu\in\sca_{\beta}$ for
$\beta=\tfrac{4}{||\mu||^2}\Ric_{\mu}$, and its eigenvalue type equals
$\beta+||\beta||^2I$ (up to a positive scalar).  In this way, the stratum to which a given
solvmanifold $S_{\mu}$ belongs determines the eigenvalue type of a potential Einstein
solvmanifold $S_{g.\mu}$ (if any).  We give some criteria to find such a stratum,
and hence the stratification provides a tool to produce existence results as well as
obstructions for nilpotent Lie algebras to be the nilradical of an Einstein
solvmanifold.

As a first application, we determine in Section \ref{dim7} all $6$-step nilpotent
Lie algebras of dimension $7$ which are Einstein nilradicals, obtaining three
explicit examples which provide a negative answer to Problem 4, as well as a new curve of examples. Secondly, we
consider in Section \ref{grafos} certain $2$-step nilpotent Lie algebras attached to
graphs and prove that they are Einstein nilradicals if and only if the graph is
positive (i.e. when certain uniquely defined weighting on the set of edges is
positive).  We next prove that any tree such that any of its edges is adjacent to at
most three other edges is positive.  This gives an existence result for Einstein
solvmanifolds. On the other hand, we prove that a graph is not positive as soon as it contains a subgraph of a certain class (see Proposition \ref{np}).  This provides a
great deal of counterexamples to Problem 4 in the $2$-step nilpotent case, the
simplest $\NN$-graded Lie algebras. Any dimension greater than or equal to $11$ is
attained.

\vs \noindent {\it Acknowledgements.}  We are very grateful to Daniel Penazzi for
helping us with graph theory and in finding examples of nonpositive graphs.  We also
thank Yuri Nikolayevsky for useful comments and suggestions.

\section{Preliminaries on Einstein solvmanifolds}\label{pre}

A {\it solvmanifold} is a simply connected solvable Lie group $S$ endowed with a
left invariant Riemannian metric.  A left invariant metric on a Lie group $G$ will
be always identified with the inner product $\ip$ determined on the Lie algebra
$\ggo$ of $G$, and the pair $(\ggo,\ip)$ will be referred to as a {\it metric Lie
algebra}.  A solvmanifold $(S,\ip)$ is called {\it standard} if $\ag:=\ngo^{\perp}$
is abelian, where $\ngo=[\sg,\sg]$, $\sg$ is the Lie algebra of $S$ and
$\ngo^{\perp}$ is relative to $\ip$.  Up to now, all known examples of noncompact
homogeneous Einstein spaces are standard Einstein solvmanifolds.  

Let us now review the results proved by Jens Heber in \cite{Hbr}.  Any solvable Lie
group admits at most one left invariant standard Einstein metric up to isometry and
scaling, and if it does, then it does not admit nonstandard Einstein metrics.  If
$(S,\ip)$ is standard Einstein, then for some distinguished element $H\in\ag$, the
eigenvalues of $\ad{H}|_{\ngo}$ are all positive integers without a common divisor,
say $k_1<...<k_r$.  If $d_1,...,d_r$ denote the corresponding multiplicities, then
the tuple
$$
(k;d)=(k_1<...<k_r;d_1,...,d_r)
$$
is called the {\it eigenvalue type} of $S$.  It turns out that $\RR H\oplus\ngo$ is
also an Einstein solvmanifold (with just the restriction of $\ip$ on it). It is then
enough to consider rank-one (i.e. $\dim{\ag}=1$) metric solvable Lie algebras as
every higher rank Einstein solvmanifold will correspond to a unique rank-one
Einstein solvmanifold and certain abelian subalgebra of derivations of $\ngo$
containing $\ad{H}$.

In every dimension, only finitely many eigenvalue types occur.  If $\mca$ is the
moduli space of all the isometry classes of Einstein solvmanifolds of a given
dimension with scalar curvature equal to $-1$, then the subspace $\mca_{{\rm st}}$
of those which are standard is open in $\mca$ (in the $C^{\infty}$-topology).  Each
eigenvalue type $(k;d)$ determines a compact path connected component $\mca_{(k;d)}$
of $\mca_{{\rm st}}$, homeomorphic to a real semialgebraic set.   

It has been proved in \cite{standard} that actually $\mca_{{\rm st}}=\mca$.
In particular, all the nice structural and uniqueness results given in \cite{Hbr}
are valid for any Einstein solvmanifold, and possibly for any noncompact homogeneous
Einstein manifold if the Alekseevskii's conjecture turns out to be true.

We now describe some interplays found in \cite{soliton, critical, praga} between
standard Einstein solvmanifolds, critical points of moment maps and Ricci soliton metrics.

Given a metric nilpotent Lie algebra $(\ngo,\ip)$, a metric solvable Lie algebra
$(\sg=\ag\oplus\ngo,\ip')$ is called a {\it metric solvable extension} of
$(\ngo,\ip)$ if the restrictions of the Lie bracket of $\sg$ and the inner product
$\ip'$ to $\ngo$ coincide with the Lie bracket of $\ngo$ and $\ip$, respectively. It
turns out that for each $(\ngo,\ip)$ there exists a unique rank-one metric solvable
extension of $(\ngo,\ip)$ which stands a chance of being an Einstein manifold (see
\cite{critical} or below). So a rank-one Einstein solvmanifold is completely
determined by its (metric) nilpotent part.  This fact turns the study of rank-one
Einstein solvmanifolds into a problem on nilpotent Lie algebras.

\begin{definition}\label{enil}
A nilpotent Lie algebra $\ngo$ is said to be an {\it Einstein nilradical} if it
admits an inner product $\ip$ such that there is a standard metric solvable
extension of $(\ngo,\ip)$ which is Einstein.
\end{definition}

In other words, Einstein nilradicals are precisely the nilradicals (i.e. the maximal
nilpotent ideal) of the Lie algebras of standard Einstein solvmanifolds.  Up to now,
the only known obstruction for a nilpotent Lie algebra to be an Einstein nilradical
is that it has to admit an $\NN$-{\it gradation}, that is, a direct sum
decomposition $\ngo=\ngo_1\oplus...\oplus\ngo_r$ (some $\ngo_i$'s might be zero)
such that $[\ngo_i,\ngo_j]\subset\ngo_{i+j}$ for all $i,j$, $i+j\leq r$ and zero
otherwise.  Such $\NN$-gradation is defined by the eigenspaces of the derivation
$\ad{H}$ mentioned above.  In this paper, we show that this is far from being the
only obstruction, by presenting several examples of $\NN$-graded (even $2$-step) Lie
algebras which are not Einstein nilradicals.

Our approach to this problem is to vary Lie brackets rather than inner products; our
main tool is the moment map for the action of the linear group on the algebraic
variety of all nilpotent Lie algebras.  This gives us the possibility to use strong
results from geometric invariant theory.

We fix an inner product vector space
$$
(\sg=\RR H\oplus\RR^n,\ip),\qquad \la H,\RR^n\ra=0,\quad \la H,H\ra=1,
$$
such that the restriction $\ip|_{\RR^n\times\RR^n}$ is the canonical inner product
on $\RR^n$, which will also be denoted by $\ip$.  A linear operator on $\RR^n$ will
be sometimes identified with its matrix in the canonical basis $\{ e_1,...,e_n\}$ of
$\RR^n$.  The metric Lie algebra corresponding to any $(n+1)$-dimensional rank-one
solvmanifold, can be modelled on $(\sg=\RR H\oplus\ngo,\ip)$ for some nilpotent Lie
bracket $\mu$ on $\RR^n$ and some $D\in\Der(\mu)$, the space of derivations of
$(\RR^n,\mu)$.  Indeed, these data define a solvable Lie bracket $[\cdot,\cdot]$ on
$\sg$ by
\begin{equation}\label{solv}
[H,X]=DX,\qquad [X,Y]=\mu(X,Y), \qquad X,Y\in\RR^n,
\end{equation}
and the solvmanifold is then the simply connected Lie group $S$ with Lie
algebra $(\sg,[\cdot,\cdot])$ endowed with the left invariant Riemannian metric
determined by $\ip$.  We shall assume from now on that $\mu\ne 0$ since the case
$\mu=0$ (i.e. abelian nilradical) is well understood (see \cite[Proposition
6.12]{Hbr}).

If $D$ is symmetric, then $(S,\ip)$ is Einstein if and only if
\begin{equation}\label{einstein}
c_{\mu}I+\tr(D)D=\Ric_{\mu},
\end{equation}
where $\Ric_{\mu}$ is the Ricci operator of $(N_{\mu},\ip)$, the simply
connected nilpotent Lie group $N_{\mu}$ with Lie algebra $(\RR^n,\mu)$ endowed with
the left invariant Riemannian metric determined by $\ip$, and
$c_{\mu}=\tfrac{\tr{\Ric_{\mu}^2}}{\tr{\Ric_{\mu}}}$ (see \cite[Lemma 2]{finding}).
Since
\begin{equation}\label{ortogonal}
\Ric_{\mu}\perp\Der(\mu)\cap\sym(n),
\end{equation}
relative to the usual inner product $\tr{\alpha\beta}$ on the space of
symmetric $n\times n $ matrices $\sym(n)$ (see \cite[(2)]{finding}), it follows from
(\ref{einstein}) that if $(S,\ip)$ is Einstein then necessarily
\begin{equation}\label{nec}
c_{\mu}I+\tr(D)D\perp\Der(\mu)\cap\sym(n).
\end{equation}

But for a given $\mu$, there exists a unique (up to a sign) symmetric derivation
$D_{\mu}$ satisfying (\ref{nec}) (possibly zero), so we can associate to each
nilpotent Lie bracket $\mu$ on $\RR^n$ a distinguished rank-one solvmanifold
$S_{\mu}:=(S_{\mu},\ip)$, defined by the data $\mu,D_{\mu}$ as in (\ref{solv}), $\tr{D_\mu}\geq 0$,
which is the only one with a chance of being Einstein among all those metric
solvable extensions of $(\mu,\ip)$.

Note that conversely, any $(n+1)$-dimensional rank-one Einstein solvmanifold is
isometric to $S_{\mu}$ for some nilpotent $\mu$ (it follows from \cite[4.10]{Hbr}
that we can assume, without loss of generality, that $\ad{H}$ is symmetric). Thus
the set $\nca$ of all nilpotent Lie brackets on $\RR^n$ parametrizes a space of
$(n+1)$-dimensional rank-one solvmanifolds
$$
\{ S_{\mu}:\mu\in\nca\},
$$
containing all those which are Einstein in that dimension.

If we consider the vector space
$$
V=\lam=\{\mu:\RR^n\times\RR^n\longrightarrow\RR^n : \mu\; \mbox{bilinear and
skew-symmetric}\},
$$
then
$$
\nca=\{\mu\in V:\mu\;\mbox{satisfies Jacobi and is nilpotent}\}
$$
is an algebraic subset of $V$ as the Jacobi identity and the nilpotency condition
can both be written as zeroes of polynomial functions.  There is a natural action of
$\G:=\Gl_n(\RR)$ on $V$ given by
\begin{equation}\label{action}
g.\mu(X,Y)=g\mu(g^{-1}X,g^{-1}Y), \qquad X,Y\in\RR^n, \quad g\in\G,\quad \mu\in V.
\end{equation}
Note that $\nca$ is $\G$-invariant and Lie algebra isomorphism classes are precisely
$\G$-orbits.  Concerning the identification $\mu\longleftrightarrow (N_{\mu},\ip)$,
this $\G$-action on $\nca$ has the following geometric interpretation: each $g\in\G$
determines a Riemannian isometry
\begin{equation}\label{id}
(N_{g.\mu},\ip)\longrightarrow (N_{\mu},\la g\cdot,g\cdot\ra)
\end{equation}
by exponentiating the Lie algebra isomorphism
$g^{-1}:(\RR^n,g.\mu)\longrightarrow(\RR^n,\mu)$.  Thus the orbit $\G.\mu$ may be
viewed as a parametrization of the set of all left invariant metrics on $N_{\mu}$.
By a result of E. Wilson, two pairs $(N_{\mu},\ip)$, $(N_{\lambda},\ip)$ are
isometric if and only if $\mu$ and $\lambda$ are in the same $\Or(n)$-orbit (see
\cite[Appendix]{minimal}), where $\Or(n)$ denotes the subgroup of $\G$ of orthogonal
matrices. Also, two solvmanifolds $S_{\mu}$ and $S_{\lambda}$ with
$\mu,\lambda\in\nca$ are isometric if and only if there exists $g\in\Or(n)$ such
that $g.\mu=\lambda$ (see \cite[Proposition 4]{critical}).  From (\ref{id}) and the
definition of $S_{\mu}$ we obtain the following result.

\begin{proposition}\label{enilmu}
If $\mu\in\nca$ then the nilpotent Lie algebra $(\RR^n,\mu)$ is an Einstein
nilradical if and only if $S_{g.\mu}$ is Einstein for some $g\in\G$.
\end{proposition}

Recall that being an Einstein nilradical is a property of a whole $\G$-orbit in
$\nca$, that is, of the isomorphism class of a given $\mu$.

The canonical inner product $\ip$ on $\RR^n$ defines an $\Or(n)$-invariant inner
product on $V$, denoted also by $\ip$, as follows:
\begin{equation}\label{innV}
\la\mu,\lambda\ra=\sum\limits_{ijk}\la\mu(e_i,e_j),e_k\ra\la\lambda(e_i,e_j),e_k\ra.
\end{equation}

\begin{theorem}\label{crit}\cite{soliton,critical}
For $\mu\in\nca$, the following conditions are equivalent:
\begin{itemize}
\item[(i)] $S_{\mu}$ is Einstein.

\item[(ii)] $\mu$ is a critical point of the functional $F:V\longrightarrow\RR$ defined by
$$
F(\mu)=\frac{16\tr{\Ric_{\mu}^2}}{||\mu||^4}.
$$

\item[(iii)] $\mu$ is a critical point of $F|_{\G.\mu}$.

\item[(iv)] $\Ric_{\mu}\in\RR I\oplus\Der(\mu)$.

\item[(v)] $\Ric_{\mu}=c_{\mu}I+D_{\mu}$.
\end{itemize}
Under these conditions, the set of critical points of $F$ lying in $\G.\mu$ equals
$\Or(n).\mu$ (up to scaling).
\end{theorem}

A Cartan decomposition for the Lie algebra $\g$ of $\G$ is given by
$\g=\sog(n)\oplus\sym(n)$, that is, in skew-symmetric and symmetric matrices
respectively.  We use the standard $\Ad(\Or(n))$-invariant inner product on $\g$,
\begin{equation}\label{inng}
\la \alpha,\beta\ra=\tr{\alpha \beta^{\mathrm t}}, \qquad \alpha,\beta\in\g.
\end{equation}

\begin{remark} We have made several abuses of notation concerning inner products.  Recall
that $\ip$ has been used to denote an inner product on $\sg$, $\RR^n$, $V$ and $\g$,
and sometimes also a left invariant metric on $S_{\mu}$ or $N_{\mu}$.
\end{remark}

The action of $\g$ on $V$ obtained by differentiation of (\ref{action}) is given by
\begin{equation}\label{actiong}
\pi(\alpha)\mu=\alpha\mu(\cdot,\cdot)-\mu(\alpha\cdot,\cdot)-\mu(\cdot,\alpha\cdot),
\qquad \alpha\in\g,\quad\mu\in V.
\end{equation}

In geometric invariant theory, a moment map for linear reductive Lie group actions
over $\CC$ has been defined in \cite{Nss} and \cite{Krw1}.  In our situation, it is
an $\Or(n)$-equivariant map
$$
m:V\smallsetminus\{ 0\}\longrightarrow\sym(n),
$$
defined implicitly by
\begin{equation}\label{defmm}
\la m(\mu),\alpha\ra=\tfrac{1}{||\mu||^2}\la\pi(\alpha)\mu,\mu\ra, \qquad \mu\in
V\smallsetminus\{ 0\}, \; \alpha\in\sym(n).
\end{equation}

Recall that $\nca\subset V$ and each $\mu\in\nca$ determines two Riemannian
manifolds $S_{\mu}$ and $(N_{\mu},\ip)$. A remarkable fact is that this moment map
encodes geometric information on $S_{\mu}$ and $(N_{\mu},\ip)$; indeed, it was
proved in \cite{minimal} that
\begin{equation}\label{mmR}
m(\mu)=\frac{4}{||\mu||^2}\Ric_{\mu}.
\end{equation}

This allows us to use strong and well-known results on the moment map due to F.
Kirwan \cite{Krw1} and L. Ness \cite{Nss}, and proved by A. Marian \cite{Mrn} in the
real case (we also refer to \cite[Section 3]{minimal} for an overview on such
results).  As a first application, we note that the functional $F$ defined in
Theorem \ref{crit}, (ii) is precisely $F(\mu)=||m(\mu)||^2$, and so we have the
following from \cite[Theorem 1, 1)]{Mrn}.

\begin{theorem}\label{minimal}
For $\mu\in\nca$, $S_{\mu}$ is Einstein if and only if $F|_{\G.\mu}$ attains its
minimum value at $\mu$.
\end{theorem}

It should be pointed out that actually most of the assertions in Theorem \ref{crit}
also follow from results proved in \cite{Mrn} (mainly from Lemma 1, Lemma 2 and
Theorem 1).

Since $\G.\mu$ parameterizes the space of left invariant metrics on $N_{\mu}$ (see
(\ref{id})) and the scalar curvature of $(N_{\mu},\ip)$ equals $-\unc||\mu||^2$, the
above result means that $S_{\mu}$ is Einstein precisely when the left invariant
metric $\ip$ on $N_{\mu}$ is very special; namely, the norm of its Ricci tensor is
minimal along all left invariant metrics on $N_{\mu}$ having the same scalar
curvature.  Such distinguished metrics are called {\it minimal} (see \cite{praga}
for further information) or sometimes {\it nilsoliton} metrics (see \cite{soliton,
Pyn}).

\begin{corollary}
A nilpotent Lie algebra is an Einstein nilradical if and only if it admits a minimal
metric.
\end{corollary}

Let $\tg$ denote the set of all diagonal $n\times n$ matrices.  If $\{
e_1',...,e_n'\}$ is the basis of $(\RR^n)^*$ dual to the canonical basis then
$$
\{ v_{ijk}=(e_i'\wedge e_j')\otimes e_k : 1\leq i<j\leq n, \; 1\leq k\leq n\}
$$
is a basis of weight vectors of $V$ for the action (\ref{action}), where $v_{ijk}$
is actually the bilinear form on $\RR^n$ defined by
$v_{ijk}(e_i,e_j)=-v_{ijk}(e_j,e_i)=e_k$ and zero otherwise.  The corresponding
weights $\alpha_{ij}^k\in\tg$, $i<j$, are given by
$$
\pi(\alpha)v_{ijk}=(a_k-a_i-a_j)v_{ijk}=\la\alpha,\alpha_{ij}^k\ra v_{ijk},
\quad\forall\alpha=\left[\begin{smallmatrix} a_1&&\\ &\ddots&\\ &&a_n
\end{smallmatrix}\right]\in\tg,
$$
where $\alpha_{ij}^k=E_{kk}-E_{ii}-E_{jj}$ and $\ip$ is the inner product defined in
(\ref{inng}).  As usual $E_{rs}$ denotes the matrix whose only nonzero coefficient
is $1$ in the entry $rs$.

For any $\mu=\sum\mu_{ij}^kv_{ijk}\in\nca$, we fix an enumeration of the set
$\left\{\alpha_{ij}^k:\mu_{ij}^k\ne 0\right\}$ and define the symmetric matrix
\begin{equation}\label{defU}
U=\left[\left\langle\alpha_{ij}^k,\alpha_{i'j'}^{k'}\right\rangle\right].
\end{equation}
We recall the following useful result.

\begin{theorem}\cite[Theorem 1]{Pyn}\label{tracy}
Assume that $\mu\in\nca$ satisfies $\Ric_{\mu}\in\tg$.  Then $S_{\mu}$ is Einstein
if and only if
$$
U\left[(\mu_{ij}^k)^2\right]=\nu [1], \qquad \nu\in\RR,
$$
where $\left[(\mu_{ij}^k)^2\right]$ is meant as a column vector in the same order
used in {\rm (}\ref{defU}{\rm )} to define $U$, and $[1]$ is the column vector
with all entries equal to $1$.
\end{theorem}

It turns out that the equations $U\left[(\mu_{ij}^k)^2\right]=\nu [1]$ are precisely
the ones given by the Lagrange method applied to find critical points of the functional $F$.

\section{A stratification for the space $\nca$}\label{st}

In this section we consider a $\G$-invariant stratification of the vector space $V$
defined in \cite{standard}. Such a stratification is an adaptation of one given by
F. Kirwan in \cite[Section 12]{Krw1} for complex reductive Lie group
representations that is strongly related to the moment map for the action.
Although the definition of the strata as well as the statement of the main theorem
are of an algebraic nature, we will show that when restricted to the space $\nca$
the stratification reveals an important interplay with geometric aspects of the
solvmanifold $S_{\mu}$ and the nilmanifold $(N_{\mu},\ip)$ attached to each
$\mu\in\nca$.  The stratum where a given solvmanifold $S_{\mu}$ lies determines the
eigenvalue type of a potential Einstein solvmanifold $S_{g.\mu}$ (if any) and so
this stratification will provide a very useful tool to verify whether a given
nilpotent Lie algebra is an Einstein nilradical or not.

Given a finite subset $X$ of $\tg$, we denote by $\CH(X)$ the convex hull of $X$ and
by $\mcc(X)$ the {\it minimal convex combination of} $X$, that is, the (unique)
vector of minimal norm (or closest to the origin) in $\CH(X)$.  Each nonzero $\mu\in
V$ uniquely determines an element $\beta_{\mu}\in\tg$ given by
$$
\beta_{\mu}=\mcc\left\{\alpha_{ij}^k:\mu_{ij}^k\ne 0\right\}, \qquad\mbox{where}\quad
\mu=\sum_{i<j}\mu_{ij}^kv_{ijk}.
$$
Recall that we always have $\beta_{\mu}\ne 0$.  Indeed, $\tr{\alpha_{ij}^k}=-1$ for
all $i<j$ and thus $\tr{\beta_{\mu}}=-1$.  Since, if $\mu$ runs through $V$, there
are only finitely many possibilities for the vectors $\beta_{\mu}$, then we can
define for each $\beta\in\tg$,
$$
\sca_{\beta}=\Big\{\mu\in V:\beta\;\mbox{is an element of maximal norm
in}\;\{\beta_{g.\mu}:g\in\G\}\Big\}.
$$
It is clear that $\sca_{\beta}$ is $\G$-invariant for any $\beta\in\tg$,
$V=\bigcup\limits_{\beta\in\tg}\sca_{\beta}$, and the set
$\{\beta\in\tg:\sca_{\beta}\ne\emptyset\}$ is finite.

For each $\beta\in\tg$ we define
$$
W_{\beta}=\{\mu\in
V:\la\beta,\alpha_{ij}^k\ra\geq||\beta||^2\quad\forall\mu_{ij}^k\ne 0\},
$$
that is, the direct sum of all the eigenspaces of $\pi(\beta)$ with eigenvalues
$\geq ||\beta||^2$.  We also consider
$$
\bca=\{\beta\in\tg^+:\sca_{\beta}\ne\emptyset\},
$$
where $\tg^+$ denotes the Weyl chamber of $\glg_n$ given by
\begin{equation}\label{weyl}
\tg^+=\left\{\left[\begin{smallmatrix} a_1&&\\ &\ddots&\\ &&a_n
\end{smallmatrix}\right]\in\tg:a_1\leq...\leq a_n\right\},
\end{equation}
and we can now state the result on the stratification.

\begin{theorem}\label{strata}\cite{standard}
The collection $\{\sca_{\beta}:\beta\in\bca\}$ is a $\G$-invariant stratification of
$V\smallsetminus\{ 0\}${\rm :}
\begin{itemize}
\item[(i)] $V\smallsetminus\{ 0\}=\bigcup\limits_{\beta\in\bca}\sca_{\beta}$ \quad {\rm (}disjoint union{\rm )}.

\item[(ii)] $\overline{\sca}_{\beta}\smallsetminus\sca_{\beta}\subset
\bigcup\limits_{||\beta'||>||\beta||}S_{\beta'}$, where $\overline{\sca}_{\beta}$ is
the closure of $\sca_{\beta}$ relative to the usual topology of $V$.  In particular,
each stratum $\sca_{\beta}$ is a locally closed subset of $V\smallsetminus\{ 0\}$.
\end{itemize}
Furthermore, for any $\beta\in\bca$ we have that
\begin{itemize}
\item[(iii)] $W_{\beta}\smallsetminus\{
0\}\subset\sca_{\beta}\cup\bigcup\limits_{||\beta'||>||\beta||}\sca_{\beta'}$.

\item[(iv)] $\sca_{\beta}\cap W_{\beta}=\{\mu\in \sca_{\beta}:\beta_{\mu}=\beta\}$.

\item[(v)] $\sca_{\beta}=\Or(n).\left(\sca_{\beta}\cap W_{\beta}\right)$.
\end{itemize}
\end{theorem}

In what follows, by using the fact that the moment map satisfies
$m(\mu)=\tfrac{4}{||\mu||^2}\Ric_{\mu}$, we derive a series of consequences of Theorem
\ref{strata} concerning the Ricci operator $\Ric_{\mu}$ of $(N_{\mu},\ip)$ and the
Einstein condition on $S_{\mu}$. We first study the behavior of $F$ relative to the
strata and give some criteria to decide, for a given $\mu\in
V$, in which stratum $\sca_{\beta}$ lies.

Since any $\mu\in\nca$ is nilpotent, the Ricci operator $\Ric_{\mu}$ of
$(N_{\mu},\ip)$ is given by (see \cite[7.38]{Bss}),
\begin{equation}\label{ricci}
\begin{array}{rl}
\la\Ric_{\mu}X,Y\ra=&-\unm\displaystyle{\sum\limits_{ij}}\la\mu(X,e_i),e_j\ra\la\mu(Y,e_i),e_j\ra \\
&+\unc\displaystyle{\sum\limits_{ij}}\la\mu(e_i,e_j),X\ra\la\mu(e_i,e_j),Y\ra,
\end{array}
\end{equation}
for all $X,Y\in\RR^n$.  We note that, in turn, the scalar curvature of
$(N_{\mu},\ip)$ is $\scalar(\mu)=\tr{\Ric_{\mu}}=-\unc||\mu||^2$.  This
formula can actually be used to define a symmetric operator $\Ric_{\mu}$ for any
$\mu\in V$.

Let $p_{\tg}(\alpha)$ denote the orthogonal projection on $\tg$ of an
$\alpha\in\sym(n)$ (i.e. the diagonal part of $\alpha$).

\begin{proposition}\label{strataF}
Let $\mu=\sum\limits_{i<j}\mu_{ij}^kv_{ijk}$ be a nonzero element of $V$.
\begin{itemize}
\item[(i)] $p_{\tg}\left(\Ri\right)=\tfrac{2}{||\mu||^2}\sum\limits_{i<j}(\mu_{ij}^k)^2
\alpha_{ij}^k\in\CH\left\{\alpha_{ij}^k:\mu_{ij}^k\ne 0\right\}$ and consequently
$$
F(\mu)\geq ||\beta_{\mu}||^2.
$$

\item[(ii)] $F(\mu)\geq ||\beta||^2$ for any $\mu\in\sca_{\beta}$.

\item[(iii)] If $\mu\in W_{\beta}$ and $\inf{F(\G.\mu)}=||\beta||^2$ then $\mu\in\sca_{\beta}$.

\item[(iv)] If $\inf{F(\G.\mu)}=||\beta_{\mu}||^2$ then $\mu\in\sca_{\beta}$, where $\beta$ is the only element
in $\tg^+$ that is conjugate to $\beta_{\mu}$.
\end{itemize}
\end{proposition}

\begin{proof} (i)  It follows from (\ref{defmm}) and (\ref{mmR}) that for any $\alpha\in\tg$,
$$
\begin{array}{rl}
\la\Ri,\alpha\ra &= \tfrac{1}{||\mu||^2}\la\pi(\alpha)\mu,\mu\ra=
\tfrac{1}{||\mu||^2}\la\sum\mu_{ij}^k\la\alpha,\alpha_{ij}^k\ra v_{ijk},\sum\mu_{ij}^kv_{ijk}\ra \\ \\
&=\tfrac{2}{||\mu||^2}\sum(\mu_{ij}^k)^2\la\alpha,\alpha_{ij}^k\ra =
\la\tfrac{2}{||\mu||^2}\sum(\mu_{ij}^k)^2\alpha_{ij}^k,\alpha\ra.
\end{array}
$$
This and the fact that $2\sum\limits_{i<j}(\mu_{ij}^k)^2=||\mu||^2$ imply the first part of (i).
For the second part it is enough to recall that
$\beta_{\mu}=\mcc\left\{\alpha_{ij}^k:\mu_{ij}^k\ne 0\right\}$ and so
$$
F(\mu)=16\tfrac{\tr{\Ric_{\mu}^2}}{||\mu||^4}=\left|\left|\Ri\right|\right|^2\geq
\left|\left| p_{\tg}\left(\Ri\right)\right|\right|^2\geq ||\beta_{\mu}||^2.
$$

\no (ii) By Theorem \ref{strata}, (iv) and (v), we can assume that
$\beta=\beta_{\mu}$ since $F$ is $\Or(n)$-invariant and hence (ii) follows directly
from (i).

\no (iii) It follows from Theorem \ref{strata}, (iii) that
$\mu\in\sca_{\beta}\cup\bigcup\limits_{||\beta'||>||\beta||}\sca_{\beta'}$ and so
from (ii) we get that necessarily $\mu\in\sca_{\beta}$ as
$\inf{F(\G.\mu)}=||\beta||^2$.

\no (iv) There exists a permutation $g\in\Or(n)$ such that $\beta_{g.\mu}=\beta$
(see the beginning of the proof of \cite[Theorem 2.10]{standard}), which implies
that $g.\mu\in W_{\beta}$ and we can apply (iii).
\end{proof}

The next step will be to describe some links with the Einstein condition given in Theorem \ref{tracy}.
Note that a linear combination $\beta_{\mu}=\sum c_{ij}^k\alpha_{ij}^k$ is not
unique in general, as the set $\left\{\alpha_{ij}^k:\mu_{ij}^k\ne 0\right\}$ might
be linearly dependent.

\begin{proposition}\label{stratatracy}
Let $\mu=\sum\limits_{i<j}\mu_{ij}^kv_{ijk}$ be a nonzero element of $V$ and consider the matrix
$U=\left[\la\alpha_{ij}^k,\alpha_{i'j'}^{k'}\ra\right]$ after fixing an enumeration
of the set $\left\{\alpha_{ij}^k:\mu_{ij}^k\ne 0\right\}$.
\begin{itemize}
    \item[(i)] If $[c_{ij}^k]$ is any solution to $U[c_{ij}^k]=\nu [1]$,
    $\nu\in\RR$, such that $\sum c_{ij}^k=1$ and all $c_{ij}^k\geq 0$, then
    $\beta_{\mu}=\sum c_{ij}^k\alpha_{ij}^k$ and $\nu=||\beta_{\mu}||^2$.

    \item[(ii)] For each convex linear combination $\beta_{\mu}=\sum c_{ij}^k\alpha_{ij}^k$, we define a finite
    set of $\lambda$'s in $V$ associated to $\mu$ by
    $\lambda=\sum\pm\sqrt{c_{ij}^k}v_{ijk}$.  If
    $\Ric_{\lambda}\in\tg$ and $\mu$ degenerates to $\lambda$ (i.e.
    $\lambda\in\overline{\G.\mu}$), then $\mu\in\sca_{\beta}$, for $\beta$ the only
    element in $\tg^+$ conjugate to $\beta_{\mu}$.

    \item[(iii)] If $\beta_{\mu}=\sum c_{ij}^k\alpha_{ij}^k$ satisfies $c_{ij}^k>0$ for any $\mu_{ij}^k\ne 0$, then
    $[c_{ij}^k]$ is a solution to $U\left[c_{ij}^k\right]=\nu [1]$ for
    $\nu=||\beta_{\mu}||^2$.
\end{itemize}
\end{proposition}

\begin{proof}
(i) It is easy to see that $U[c_{ij}^k]=\nu [1]$ is precisely the linear system from
the Lagrange multiplier method applied to find critical points of the functional
$[c_{ij}^k]\mapsto||\alpha||^2$, $\alpha=\sum c_{ij}^k\alpha_{ij}^k$, restricted to
the leaf $\sum c_{ij}^k=1$.  Since the set $\{\alpha=\sum c_{ij}^k\alpha_{ij}^k:\sum
c_{ij}^k=1\}$ is a linear variety there is a unique critical point $\alpha_m$ which
is a global minimum, and so $\sum c_{ij}^k\alpha_{ij}^k=\alpha_m$ for any solution
$[c_{ij}^k]$. If in addition all $c_{ij}^k\geq 0$ then $\sum
c_{ij}^k\alpha_{ij}^k\in\CH\left\{\alpha_{ij}^k:\mu_{ij}^k\ne 0\right\}$ and hence
$\alpha_m$ has to be $\beta_{\mu}$.  Moreover,
$$
||\beta_{\mu}||^2=\left|\left|\sum
c_{ij}^k\alpha_{ij}^k\right|\right|^2=\left\langle
U[c_{ij}^k],[c_{ij}^k]\right\rangle=\nu\sum c_{ij}^k=\nu.
$$

\no (ii) It follows from Proposition \ref{strataF}, (i) that for any such $\lambda$,
$$
\tfrac{4}{||\lambda||^2}\Ric_{\lambda}=
p_{\tg}\left(\tfrac{4}{||\lambda||^2}\Ric_{\lambda}\right)
=\tfrac{2}{||\lambda||^2}\sum c_{ij}^k\alpha_{ij}^k=\beta_{\mu},
$$
(recall that $||\lambda||^2=2$) and thus $\inf{F(\G.\mu)}\leq
F(\lambda)=||\beta_{\mu}||^2$.  Now Proposition \ref{strataF}, (iii) implies that
$\mu\in\sca_{\beta}$.

\no (iii) $\beta_{\mu}$ is the closest point to the origin in
$\CH\left\{\alpha_{ij}^k:\mu_{ij}^k\ne 0\right\}$, so if all the $c_{ij}^k$'s are
positive then $\beta_{\mu}$ is actually a local minimum of the functional mentioned
in the proof of (i).  This implies that $[c_{ij}^k]$ is a critical point of this
functional and consequently a solution to $U\left[c_{ij}^k\right]=\nu [1]$.
\end{proof}

We finally show the interplay between the stratification and the Einstein condition
on a solvmanifold $S_{\mu}$.

\begin{proposition}\label{strataEinstein}
Let $\mu\in\nca$, $\mu\ne 0$.
\begin{itemize}
\item[(i)] If $S_{\mu}$ is Einstein then $\mu\in\sca_{\beta}$ for $\beta$ the only element in $\tg^+$ conjugate
to $\Ri$.  In such case, the eigenvalue type of $S_{\mu}$ is a positive scalar
multiple of $\beta+||\beta||^2I$.

\item[(ii)] For $\mu\in\sca_{\beta}$ the following conditions are equivalent:
\begin{itemize}
\item[(a)] $S_{\mu}$ is Einstein.

\item[(b)] $\Ri$ is conjugate to $\beta$.

\item[(c)] $F(\mu)=||\beta||^2$.
\end{itemize}
\end{itemize}
\end{proposition}

\begin{proof} (i) There exists $g\in\Or(n)$ such that
$\tfrac{4\Ric_{g.\mu}}{||g.\mu||^2}=g\tfrac{4\Ric_{\mu}}{||\mu||^2}g^{-1}=\beta$.
Thus $||\beta||^2=-\tfrac{4c_{\mu}}{||g.\mu||^2}$ and it then follows from Theorem
\ref{crit}, (v) that $\beta+||\beta||^2I\in\Der(g.\mu)$.  By using that
$\la\pi(\beta+||\beta||^2I)g.\mu,g.\mu\ra=0$, it is easy to check that
$\la\beta,\alpha_{ij}^k\ra=||\beta||^2$ for any $\mu_{ij}^k\ne 0$, that is,
$g.\mu\in W_{\beta}$. Thus $g.\mu\in\sca_{\beta}$ (and so $\mu\in\sca_{\beta}$) by
Proposition \ref{strataF}, (iii) since $F(g.\mu)=||\beta||^2$.

\no (ii) It follows from (i) that (a) implies (b), and (c) follows from (b)
trivially.  If we assume (c) then Proposition \ref{strataF}, (ii) implies that
$F(g.\mu)\geq||\beta||^2=F(\mu)$ for any $g\in\G$ since $\sca_{\beta}$ is
$\G$-invariant, and hence $S_{\mu}$ is Einstein by Theorem \ref{minimal}.
\end{proof}

In the light of Theorem \ref{crit}, another natural approach to find rank-one
Einstein solvmanifolds would be to use the negative gradient flow of the functional
$F$.  It follows from \cite[Lemma 6]{critical} that
$$
\grad(F)_{\mu}=-\tfrac{16}{||\mu||^6}\left(\delta_{\mu}(\Ric_{\mu})||\mu||^2+4\tr{\Ric_{\mu}^2}\mu\right),
$$
where $\delta_{\mu}:\g\longrightarrow V$ is defined by
$\delta_{\mu}(\alpha)=-\pi(\alpha)\mu$ (see (\ref{actiong})).  Since $F$ is
invariant under scaling we know that $||\mu||$ will remain constant in time along
the flow. We may therefore restrict ourselves to the sphere of radius $2$, where the
negative gradient flow $\mu=\mu(t)$ of $F$ becomes
\begin{equation}\label{flow}
    \ddt\mu=\delta_{\mu}(\Ric_{\mu})+\tr{\Ric_{\mu}^2}\mu.
\end{equation}

Notice that $\mu(t)$ is a solution to this ODE if and only if
$g.\mu(t)$ is so for any $g\in\Or(n)$ (use that
$g.\delta_{\mu}(\Ric_{\mu})=\delta_{g.\mu}(g\Ric_{\mu}g^{-1})=\delta_{g.\mu}(\Ric_{g.\mu})$),
in accordance with the $\Or(n)$-invariance of $F$.  The existence of a solution $\mu(t)$, $t\in[0,\infty)$, is guaranteed by the compactness of the sphere, and the existence of a unique limit $\lim\limits_{t\to\infty}\mu(t)$ follows from the fact that $F$ is a polynomial (see \cite{Krd} or \cite{Mss}).

\begin{proposition}\label{strataflow}
For $\mu_0\in V$, $||\mu_0||=2$, let $\mu(t)$ be the flow defined in {\rm
(\ref{flow})} with $\mu(0)=\mu_0$ and put $\lambda=\lim\limits_{t\to\infty}\mu(t)$.
Then
\begin{itemize}
    \item[(i)] $S_{\lambda}$ is Einstein.
    \item[(ii)] $\lambda\in\overline{\G.\mu_0}$ (i.e. $\mu_0$ degenerates to $\lambda$).
    \item[(iii)] If $\beta:=\Ric_{\lambda}\in\tg^+$ and $\mu_0\in W_{\beta}$ then $\mu_0\in\sca_{\beta}$.
\end{itemize}
\end{proposition}

\begin{proof}
Part (i) follows from Theorem \ref{crit} by using that $\lambda$ is a critical point of
$F$. Since $\grad(F)_{\mu}\in\tang_{\mu}\G.\mu$ for any $\mu\in V$ we have that
$\mu(t)\in\G.\mu_0$ for all $t$, which proves (ii).  For (iii), we just apply (ii)
and Proposition \ref{strataF}, (iii).
\end{proof}

In order to show the interplay between the geometry of the moment map and algebra, we will apply the results of this section in two specific cases.  More involved applications will be given in Sections
\ref{dim7} and \ref{grafos}.

\begin{example}\label{flow1}
Let $(\RR^{7},\mu)$ be the $7$-dimensional $2$-step nilpotent Lie algebra defined by
$$
\begin{array}{lll}
  \mu_0(e_1,e_2)=\sqrt{\tfrac{2}{3}}e_5, & \mu_0(e_2,e_3)=\sqrt{\tfrac{2}{3}}e_6, &
  \mu_0(e_3,e_4)=\sqrt{\tfrac{2}{3}}e_7.
\end{array}
$$
The scaling by $\sqrt{\tfrac{2}{3}}$ is just to get $||\mu_0||=2$.  A curve
$\mu=\mu(t)$ of the form
$$
\begin{array}{lll}
  \mu(e_1,e_2)=a(t)e_5, & \mu(e_2,e_3)=b(t)e_6, &
  \mu(e_3,e_4)=c(t)e_7,
\end{array}
$$
satisfies (\ref{flow}) if and only if
$$
\begin{array}{l}
  a'=-\unm a(3a^2+b^2)+\unc(3a^4+3b^4+3c^4+2a^2b^2+2b^2c^2)a, \\ \\
  b'=-\unm b(a^2+3b^2+c^2)+\unc(3a^4+3b^4+3c^4+2a^2b^2+2b^2c^2)b, \\ \\
  c'=-\unm c(b^2+3c^2)+\unc(3a^4+3b^4+3c^4+2a^2b^2+2b^2c^2)c. \\
\end{array}
$$
Assume that $\mu(0)=\mu_0$, that is, $a(0)=b(0)=c(0)=\sqrt{\tfrac{2}{3}}$.  By
subtracting the last equation to the first one we obtain that $a(t)=c(t)$ for all
$t$, and then by letting $x:=a^2$, $y=b^2$, we get the equivalent system
$$
\begin{array}{l}
  x'=x\left(\unm(6x^2+3y^2+4xy)-3x-y\right), \qquad x(0)=\tfrac{2}{3}, \\ \\
  y'=y\left(\unm(6x^2+3y^2+4xy)-2x-3y\right), \qquad y(0)=\tfrac{2}{3}. \\
\end{array}
$$
This implies that
$$
\left(\frac{x}{y}\right)'=\frac{x}{y}(2y-x), \qquad \frac{x}{y}(0)=1, \quad
\left(\frac{x}{y}\right)'(0)=\tfrac{4}{3}>0,
$$
and hence $\lim\limits_{t\to\infty}x(t)=2\lim\limits_{t\to\infty}y(t)$.  Since
$4x+2y=||\mu||^2=4$ for all $t$ we conclude that
$\lambda=\lim\limits_{t\to\infty}\mu(t)$ is given by
$$
\begin{array}{lll}
  \lambda(e_1,e_2)=\tfrac{2}{\sqrt{5}}e_5, & \lambda(e_2,e_3)=\sqrt{\tfrac{2}{5}}e_6, &
  \lambda(e_3,e_4)=\tfrac{2}{\sqrt{5}}e_7.
\end{array}
$$
Evidently $\lambda\in\Gl_7.\mu_0$ and so $\mu_0$ is an Einstein nilradical by
Proposition \ref{strataflow}, (i) and Proposition \ref{enilmu}.
\end{example}

We now use the geometric technique to present an example of a $2$-step nilpotent Lie algebra which is
not an Einstein nilradical.  This is the first known example of this kind to our best knowledge.

\begin{example}\label{flow2}
Let $(\RR^{11},\mu_0)$ be the $11$-dimensional $2$-step nilpotent Lie algebra
defined by
$$
\begin{array}{lll}
  \mu_0(e_1,e_2)=\sqrt{\tfrac{2}{5}}e_7, & \mu_0(e_1,e_4)=\sqrt{\tfrac{2}{5}}e_9, & \mu_0(e_2,e_6)=\sqrt{\tfrac{2}{5}}e_{11}. \\
  \mu_0(e_1,e_3)=\sqrt{\tfrac{2}{5}}e_8, & \mu_0(e_2,e_5)=\sqrt{\tfrac{2}{5}}e_{10}, &  \\
\end{array}
$$
A curve $\mu=\mu(t)$ of the form
$$
\begin{array}{lll}
  \mu_0(e_1,e_2)=b(t)e_7, & \mu_0(e_1,e_4)=a(t)e_9, & \mu_0(e_2,e_6)=a(t)e_{11}, \\
  \mu_0(e_1,e_3)=a(t)e_8, & \mu_0(e_2,e_5)=a(t)e_{10}, &  \\
\end{array}
$$
satisfies (\ref{flow}) if and only if
$$
\begin{array}{l}
  a'=-\unm a(4a^2+b^2)+\unc(16a^4+3b^4+8a^2b^2)a, \\ \\
  b'=-\unm b(4a^2+3b^2)+\unc(16a^4+3b^4+8a^2b^2)b. \\
\end{array}
$$
Assume that $\mu(0)=\mu_0$, that is, $a(0)=b(0)=\sqrt{\tfrac{2}{5}}$ and let
$x:=a^2$, $y=b^2$, to get the equivalent system
$$
\begin{array}{l}
  x'=x\left(\unm(16x^2+3y^2+8xy)-4x-y\right), \qquad x(0)=\tfrac{2}{5}, \\ \\
  y'=y\left(\unm(16x^2+3y^2+8xy)-4x-3y\right), \qquad y(0)=\tfrac{2}{5}. \\
\end{array}
$$
This implies that
$$
\left(\frac{y}{x}\right)'=\frac{y}{x}(-2y), \qquad \frac{y}{x}(0)=1, \quad
\left(\frac{y}{x}\right)'(0)=-\tfrac{4}{5}<0,
$$
and hence $\lim\limits_{t\to\infty}y(t)=0$.  Since $8x+2y=||\mu||^2=4$ for all $t$
we obtain that $\lim\limits_{t\to\infty}x(t)=\unm$ and
$\lambda=\lim\limits_{t\to\infty}\mu(t)$ is therefore given by
$$
\begin{array}{lll}
  \lambda(e_1,e_2)=0, & \lambda(e_1,e_4)=\tfrac{1}{\sqrt{2}}e_9, & \lambda(e_2,e_6)=\tfrac{1}{\sqrt{2}}e_{11}, \\
  \lambda(e_1,e_3)=\tfrac{1}{\sqrt{2}}e_8, & \lambda(e_2,e_5)=\tfrac{1}{\sqrt{2}}e_{10}, &  \\
\end{array}
$$
Thus $\lambda\notin\Gl_{11}.\mu_0$ since its derived algebra is $4$-dimensional, one
less than for $\mu_0$.  An easy computation yields that $\beta:=\Ric_{\lambda}$ is
the element in $\tg^+$ with entries
$$
-\unm,-\unm,-\unc,-\unc,-\unc,-\unc,0,\unc,\unc,\unc,\unc.
$$
It is not hard to check that $\mu_0\in W_{\beta}$ and so $\mu_0\in\sca_{\beta}$ by
Proposition \ref{strataflow}, (iii).  Assume now that $\mu_0$ is an Einstein
nilradical.  Then there must exist $\mu\in\Gl_{11}.\mu_0$ such that $S_{\mu}$ is
Einstein (see Proposition \ref{enilmu}) and so $\frac{4}{||\mu||^2}\Ric_{\mu}$ is
necessarily conjugate to $\beta$ by Proposition \ref{strataEinstein}, (ii).  This
implies that $\Ric_{\mu}$ has zero as an eigenvalue, a contradiction by Lemma
\ref{2step} since the derived algebra of $\mu$ must coincide with its center (recall
that $\mu\simeq\mu_0$).
\end{example}

To conclude this section, we give some general results on the Ricci operator of
nilmanifolds which will be useful in the results in the next sections.

\begin{lemma}\label{2step}
Let $(\ngo,\ip)$ be a metric $2$-step nilpotent Lie algebra.  Then $[\ngo,\ngo]$
coincides with the center of $\ngo$ if and only if the eigenvalues of the Ricci
operator $\Ric_{\ip}$ are all nonzero.
\end{lemma}

\begin{proof}
Consider the orthogonal decomposition of the form $\ngo=\vg\oplus[\ngo,\ngo]$.  It
follows from (\ref{ricci}) that $\Ric_{\ip}$ leaves $\vg$ and $[\ngo,\ngo]$
invariant and it is positive definite on $[\ngo,\ngo]$ and negative semidefinite on
$\vg$.  Moreover, $\Ric_{\ip}X=0$ for $X\in\vg$ if and only if $X$ is in the center
of $\ngo$.  This concludes the proof.
\end{proof}

\begin{lemma}\label{fisica}
For $\mu\in V$, the following two conditions are sufficient to have $\Ric_{\ip}\in\tg$:
\begin{itemize}
    \item for all $i<j$ there is at most one $k$ such that $\mu_{ij}^k\ne 0$,
    \item if $\mu_{ij}^k$ and $\mu_{i'j'}^k$ are nonzero then either $\{ i,j\}=\{
    i',j'\}$ or $\{ i,j\}\cap\{ i',j'\}=\emptyset$.
\end{itemize}
\end{lemma}

\begin{proof}
The lemma follows directly from (\ref{ricci}).
\end{proof}

\section{Applications in dimension $7$}\label{dim7}

It is known that, up to dimension $6$, any nilpotent Lie algebra is an Einstein
nilradical (see \cite{finding,Wll}), thus $7$ is the first dimension to consider Problem 4.  We will
determine in this section the $7$-dimensional $6$-step Einstein
nilradicals, obtaining in particular three examples which can not be so.

If $\mu\in V=\Lambda^2(\RR^7)^*\otimes\RR^7$ satisfies $\mu(e_i,e_j)=a_{ij}e_{i+j}$,
$i<j$, and zero otherwise, then it will be denoted by the $9$-tuple
\begin{equation}\label{gradedmu}
\mu=(a_{12}, a_{13},a_{14},a_{15},a_{16},a_{23},a_{24},a_{25},a_{34})
\end{equation}
(i.e. $\mu_{ij}^{i+j}=a_{ij}$).  Assume that all the $a_{ij}$'s are nonzero.  Thus
the set $\left\{\alpha_{ij}^k:\mu_{ij}^k\ne 0\right\}$ is given by
$$
\{ (-1,-1,1,0,0,0,0), (-1,0,-1,1,0,0,0), ..., (0,0,-1,-1,0,0,1)\},
$$
and with respect to this enumeration the matrix $U$ defined in (\ref{defU}) is
$$
U=\left[%
\begin{array}{ccccccccc}
3 & 0 & 1 & 1 & 1 & 0 & 1 & 1 & -1 \\ 0 & 3 & 0 & 1 & 1 & 1 & -1 & 0 & 0 \\ 1 & 0 & 3 & 0 & 1 & 1 & 1 & -1 & 1 \\
1 & 1 & 0 & 3 & 0 & -1 & 1 & 1 & 0 \\ 1 & 1 & 1 & 0 & 3 & 0 & -1 & 1 & 1 \\
0 & 1 & 1 & -1 & 0 & 3 & 1 & 0 & 1 \\ 1 & -1 & 1 & 1 & -1 & 1 & 3 & 1 & 1 \\ 1 & 0 & -1 & 1 & 1 & 0 & 1 & 3 & 1 \\
-1 & 0 & 1 & 0 & 1 & 1 & 1 & 1 & 3
\end{array}%
\right].
$$
It is easy to check that $\Ric_{\mu}\in\tg$ for any $\mu$ of the form
(\ref{gradedmu}) by using Lemma \ref{fisica}.  It follows from Theorem \ref{tracy}
that $S_{\mu}$ is Einstein if
and only if $U\left[\begin{smallmatrix} a_{12}^2 \\ \vdots \\
a_{34}^2\end{smallmatrix}\right]=\nu\left[\begin{smallmatrix} 1 \\ \vdots \\
1\end{smallmatrix}\right]$ for some $\nu\in\RR$.  We solve this linear system and
get that $S_{\mu}$ is Einstein if and only if there exists $a,b,c\in\RR$ such that
$$
\begin{smallmatrix}
\mu=\mu_{a,b,c} := \left(\pm\sqrt{a}, \pm\sqrt{2-b}, \pm\sqrt{3-a-b-c}, \pm\sqrt{b},
\pm\sqrt{b+c-1}, \pm\sqrt{b}, \pm\sqrt{c}, \pm\sqrt{3-a-b-c}, \pm\sqrt{a}\right)
\end{smallmatrix}
$$
up to a scalar multiple, where all the numbers under a square root must be of course
nonnegative.  Recall that these are precisely the critical points of the functional
$F(a_{12},...,a_{34})=\tfrac{16\tr{\Ric_{\mu}^2}}{||\mu||^2}$ restricted to any leaf
of the form $\sum a_{ij}^2=$ constant.  By Proposition \ref{stratatracy}, (i) we
have that for any $\mu=\mu_{a,b,c}$,
$$
\beta_{\mu}=\tfrac{1}{7}\left(a\alpha_{12}^3+(2-b)\alpha_{13}^4+...+a\alpha_{34}^7\right)
=\tfrac{1}{7}(-4,-3,-2,-1,0,1,2).
$$
This element of $\tg^+$ will be denoted by $\beta$ from now on.  Thus $\Ri=\beta$
(see Proposition \ref{strataF}, (i)), $\mu\in\sca_{\beta}$ (see Proposition
\ref{strataEinstein}, (i)) and the corresponding eigenvalue type equals
$(1<2<3<4<5<6<7;1,...,1)$ for any $\mu=\mu_{a,b,c}$.  Recall that the Jacobi
identity adds some conditions on $a,b,c$.

It follows from \cite[Theoreme 2]{GzHkm} and \cite[Theorem 5.17]{Mll} that any
$7$-dimensional $6$-step nilpotent real Lie algebra admitting an $\NN$-gradation is
isomorphic to one and only one of the following:
\begin{itemize}
    \item $\mu_1=(1,...,1,0,0,0,0)$,
    \item[ ]
    \item $\mu_2=(1,...,1,0,1,1,1)$,
    \item[ ]
    \item $\lambda_t=(1,...,1,t,1-t), \qquad t\in\RR$,
    \item[]
    \item $\mu_3=\mu_1+v_{237}$,
    \item[]
    \item $\mu_4=\mu_1+v_{236}+v_{247}$.
\end{itemize}
It should be noticed that this list differs a little from the one given in
\cite{Mll}.  The translation is as follows: $\mu_1=m_0(7)$; $\mu_2\simeq
m_{0,1}(7)=\ggo_{7,-2}$; $\lambda_t\simeq\ggo_{7,\alpha}$,
$\alpha=\tfrac{2t-1}{1-t}$, $t\ne 1$; $\lambda_1=m_2(7)$.  We also recall that
although the ground field in \cite{GzHkm} is $\CC$, the result we need is still
valid on $\RR$ since it only uses Lemme 2, which is clearly valid for the real
numbers as well.

Let $g\in\Gl_7$ be the diagonal matrix with entries
$$
\begin{array}{lcl}
  g_1=1, && g_5=\tfrac{\sqrt{a}(3-a-b-c)(2-b)}{\sqrt{b}}, \\ \\
  g_2=\tfrac{\sqrt{3-a-b-c}\sqrt{2-b}}{\sqrt{b}}, &&
  g_6=\sqrt{a}(3-a-b-c)(2-b), \\ \\
  g_3=\tfrac{\sqrt{a}\sqrt{3-a-b-c}\sqrt{2-b}}{\sqrt{b}}, &&
  g_7=\tfrac{\sqrt{a}(3-a-b-c)\sqrt{2-b}\sqrt{b+c-1}}{\sqrt{c}}. \\ \\
  g_4=\tfrac{\sqrt{a}\sqrt{3-a-b-c}(2-b)}{\sqrt{b}}, &&
\end{array}
$$
It is a straightforward calculation to prove that if $c=\tfrac{b^2}{2-b}$ then
$$
g^{-1}.\mu_{a,b,c}=\left(1,...,1,\pm\tfrac{b(3-a-b-c)}{c\sqrt{b+c-1}\sqrt{2-b}},
\pm\tfrac{a}{\sqrt{c}\sqrt{b+c-1}}\right),
$$
where a positive square root has already been chosen in all the coordinates of
$\mu_{a,b,c}$ except in the last two.  It is not hard to check that
$g^{-1}.\mu_{a,b,c}$ is as follows depending on $a,b,c$:
\begin{itemize}
    \item $\lambda_t$, $0<t<1$, if $t=1-a$, $b=c=1$ \quad (positive square roots chosen in $\mu_{a,b,c}$);
    \item[ ]
    \item $\lambda_t$, $1<t<\infty$, if $t=\tfrac{b(3-a-b-c)}{c\sqrt{b+c-1}\sqrt{2-b}}$,
    $1-t=-\tfrac{a}{\sqrt{c}\sqrt{b+c-1}}$, $a=\tfrac{6-5b-b\sqrt{3b-2}}{2(2-b)}$,
    $c=\tfrac{b^2}{2-b}$, $\tfrac{2}{3}<b<1$;
    \item[ ]
    \item $\lambda_t$, $-\infty<t<0$, if $t=-\tfrac{b(3-a-b-c)}{c\sqrt{b+c-1}\sqrt{2-b}}$,
    $1-t=\tfrac{a}{\sqrt{c}\sqrt{b+c-1}}$, $a=\tfrac{6-5b+b\sqrt{3b-2}}{2(2-b)}$,
    $c=\tfrac{b^2}{2-b}$, $\tfrac{2}{3}<b<1$.
\end{itemize}
This implies that $\lambda_t$ is an Einstein nilradical for any $t\ne 0,1$.  On the
other hand, it is proved in \cite[Theorem 4.2]{finding} that $\mu_1$ is an Einstein
nilradical, as it is isomorphic to
$\mu_1'=(\sqrt{5},\sqrt{8},3,\sqrt{8},\sqrt{5},0,0,0,0)$ with $S_{\mu_1'}$ Einstein
of eigenvalue type $(1<16<17<18<19<20<21;1,...,1)$ (see also \cite[Theorem
27]{Pyn}). We also have that $\mu_3$ and $\mu_4$ are both Einstein nilradicals as
they are isomorphic to the critical points
$$
\begin{array}{l}
\mu_3'=(\sqrt{5},\sqrt{5},3,\sqrt{8},\sqrt{2},0,0,0,0)+3v_{237}, \\ \\
\mu_4'=(\sqrt{10},\sqrt{21},\sqrt{18},4,\sqrt{18},0,0,0,0)+\sqrt{21}v_{236}+\sqrt{18}v_{247},
\end{array}
$$
of eigenvalue type $(1<4<5<6<7<8<9;1,...,1)$ and $(1<3<4<5<6<7<8;1,...,1)$,
respectively. Recall that $\mu_1',\mu_3',\mu_4'$ are all different from any
$\mu_{a,b,c}$, in accordance with the fact that they belong to other strata.

It then only remains to consider the cases $\mu_2$, $\lambda_0$ and $\lambda_1$. Let
us first consider $\mu_2$. If $U_1$ is the matrix associated to the ordered set of
weights $\{\alpha_{12}^3,...,\hat{\alpha_{23}^5},\alpha_{34}^7\}$ (i.e. only
$\alpha_{23}^5$ is missing), then the solutions to $U_1[c_{ij}^k]=\nu [1]$ are given
up to a scalar multiple by
$$
(a,2,3-a-c,0,c-1,c,3-a-c,a).
$$
It follows from Proposition \ref{stratatracy}, (i) that
$\beta_{\mu_2}=\beta=\tfrac{1}{7}(-4,...,2)$ and hence in order to apply part (ii)
of the same proposition we should try to find a degeneration
$\mu_2\rightarrow\lambda$ for some of the $\lambda's$ defined there, that is,
$\lambda=\mu_{a,0,c}$ for some $a,c\in\RR$ such that $a,3-a-c,c\geq 0$.  This may be
done by defining a curve $g_t\in\Gl_7$ of diagonal matrices such that
$\lim\limits_{t\rightarrow 0} g_t.\mu_2=\mu_{a,0,c}$.  By assuming that the square
of the entries of $g_t^{-1}$ are of the form $s_it^{r_i}$, we get the following
necessary and sufficient conditions to get such a degeneration:
$$
\begin{array}{lll}
\frac{s_1s_2}{s_3}=a, & \frac{s_2s_4}{s_6}=c, & r_i+r_j=r_{i+j}\quad\forall(i,j)\ne (1,5),(2,3), \\ \\
\frac{s_1s_3}{s_4}=2, & \frac{s_2s_5}{s_7}=3-a-c, & r_1+r_5>r_6. \\ \\
\frac{s_1s_4}{s_5}=3-a-c, & \frac{s_3s_4}{s_7}=a, & \\ \\ \frac{s_1s_6}{s_7}=c-1, &
&
\end{array}
$$
It is easy to see that this is equivalent to $c=\tfrac{\sqrt{28}-1}{3}$,
$a=\sqrt{c(c-1)}$, all the $r_i$'s are in terms of $r_1,r_2$ and $2r_1>r_2$.  So
that we can take $g_t$ with entries
$$
1,t,\sqrt{a}t,\sqrt{2a}t,\sqrt{2a(3-a-c)}t,\sqrt{2ac}t^2,\sqrt{2a^3}t^2,
$$
with $a,c$ as above.  It then follows from Proposition \ref{stratatracy}, (ii) that
$\mu_2\in\sca_{\beta}$.

We now suppose that $\mu_2$ is an Einstein nilradical, that is, there exists
$\mu\in\Gl_7.\mu_2$ such that $S_{\mu}$ is Einstein (see Proposition \ref{enilmu}).
By Proposition \ref{strataEinstein} we can assume that
$$
\beta+||\beta||^2I=\tfrac{1}{7}\left[
\begin{smallmatrix}
  1 &  &  &  &  &  &  \\
   & 2 &  &  &  &   &   \\
    &   & 3 &   &   &   &   \\
    &   &   & 4 &   &   &   \\
    &   &   &   & 5 &   &   \\
    &   &   &   &   & 6 &   \\
    &   &   &   &   &   & 7
\end{smallmatrix}
\right]
$$
is a derivation of $\mu$.  This implies that $\mu$ has the form (\ref{gradedmu}) and
since $S_{\mu}$ is Einstein there must exist $a,b,c$ such that $\mu=\mu_{a,b,c}$.
Recall that $\mu$ is isomorphic to $\mu_2$.  Thus $\mu$ is $6$-step nilpotent, which
implies that $\mu_{1j}^{1+j}\ne 0$ for $j=2,...,6$ and so $\mu_{23}^5$ is nonzero as
well since $\mu_{23}^5=\pm\sqrt{b}=\pm\mu_{15}^6$.  We finally arrive at a
contradiction since the Lie algebra invariant given by $\ggo/C_4(\ggo)$, where
$C_4(\ggo)$ is the fourth term in the descending central series of the Lie algebra
$\ggo$, gives rise to non-isomorphic $5$-dimensional Lie algebras when computed for
$\mu$ and $\mu_2$, respectively (note that $\mu_{23}^5\ne 0=(\mu_2)_{23}^5$ and see
\cite{Mgn}).  We then conclude that $\mu_2$ is not an Einstein nilradical.

Since $\lambda_0\in W_{\beta}$ we can prove that $\lambda_0\in\sca_{\beta}$ by using
the following degeneration: $\lim\limits_{t\rightarrow\infty}
g_t.\lambda_0=\mu_{1,1,1}$, $g_t=e^{-t\alpha}$, $\alpha\in\tg$ with entries
$(0,1,1,1,2,2,2)$.  Indeed, this shows that $\inf{F(\Gl_7.\lambda_0)}\leq
||\beta||^2=\tfrac{5}{7}=F(\mu_{1,1,1})$ and we can apply Proposition \ref{strataF},
(iii).  Once we know that $\lambda_0\in\sca_{\beta}$ then we argue as in the case of
$\mu_2$ above and get that $\lambda_0$ must be isomorphic to a $\mu=\mu_{a,b,c}$ if
we assume that $\lambda_0$ is an Einstein nilradical.  We use this time the
invariant $\dim\{\ad{X}|_{C_3(\ggo)}:X\in\ggo\}$, which equals $1$ for $\lambda_0$
and equals $2$ for $\mu$ since $(\lambda_0)_{25}^7=0$ and
$\mu_{25}^7=\pm\mu_{14}^5\ne 0$ ($\mu$ $6$-step), respectively.

Finally, for $\lambda_1$ we argue in an analogous way, by using the degeneration
$\lim\limits_{t\rightarrow\infty} g_t.\lambda_1=\mu_{0,1,1}$, $g_t=e^{-t\alpha}$,
$\alpha\in\tg$ with entries $(1,2,0,1,2,3,4)$, and the invariant $C_1(\ggo)$, which
is abelian for $\lambda_1$ and nonabelian for any $\mu_{a,b,c}$ with $a\ne 0$.

The results obtained in this section may be summarized in the following
classification.

\begin{theorem}\label{clasi7}
A $6$-step nilpotent Lie algebra $\mu$ of dimension $7$ is an Einstein nilradical if and
only if $\mu\ne \mu_2, \lambda_0, \lambda_1$.
\end{theorem}

\section{Application to $2$-step nilpotent Lie algebras attached to
graphs}\label{grafos}

Let $\gca$ be a graph with set of vertices $\{ v_1,...,v_p\}$ and set of edges $\{
l_1,...,l_q\}\subset\{ v_iv_j:1\leq i,j\leq p\}$.  We associate to $\gca$ a $2$-step
nilpotent Lie algebra $(\RR^n,\mu_{\gca})$, $n=p+q$, defined by
$$
\mu_{\gca}(e_i,e_j)=
\left\{%
\begin{array}{ll}
    e_{p+k}, & \hbox{if $l_k=v_iv_j$, $i<j$;} \\ \\
    0, & \hbox{{\rm otherwise}.} \\
\end{array}%
\right.
$$
Thus the center of $\mu_{\gca}$ coincides with the derived algebra
$\mu_{\gca}(\RR^n,\RR^n)$, which is linearly generated by $e_{p+1},...,e_n$, if and
only if $\gca$ has no any isolated point.  These Lie algebras have been recently
considered in \cite{DnMnk}, where their automorphism groups have been computed.  Our
aim in this section is to consider the question of for which graphs $\gca$ the Lie
algebra $\mu_{\gca}$ is an Einstein nilradical.

For $a_1,...,a_q\in\RR$ we consider the Lie algebra $\mu_{\gca}(a_1,...,a_q)$
defined by
$$
\mu_{\gca}(a_1,...,a_q)(e_i,e_j)=\left\{%
\begin{array}{ll}
    a_ke_{p+k}, & \hbox{if $l_k=v_iv_j$, $i<j$;} \\ \\
    0, & \hbox{{\rm otherwise}.} \\
\end{array}%
\right.
$$
Thus $\mu_{\gca}=\mu_{\gca}(1,...,1)$ and $\mu_{\gca}$ is isomorphic to
$\mu_{\gca}(a_1,...,a_k)$ if and only if $a_k\ne 0$ for any $k=1,...,q$.  Recall
that $\left(\mu_{\gca}\right)_{ij}^{p+k}=1$ if $l_k=v_iv_j$ and zero otherwise. Thus
$\beta_{\mu_{\gca}}$ can be written as
$$
\beta_{\mu_{\gca}}=\sum_{l_k=v_iv_j}c_k\alpha_{ij}^{p+k}, \qquad c_1,...,c_q\geq 0,
\quad \sum_{k=1}^qc_k=1.
$$
If $\lambda=\mu_{\gca}(\sqrt{c_1},...,\sqrt{c_q})$ then $\mu_{\gca}$ degenerates to
$\lambda$ since
$\lim\limits_{t\to\infty}\mu_{\gca}(\tfrac{1}{t}+\sqrt{c_1},...,\tfrac{1}{t}+\sqrt{c_1})=\lambda$
and $\tfrac{1}{t}+\sqrt{c_i}>0$ for all $i$ and $t>0$.  We may therefore apply
Proposition \ref{stratatracy}, (ii) since $\Ric_{\lambda}\in\tg$ by Lemma
\ref{fisica} and obtain:

\begin{proposition}\label{gt}
For any graph $\gca$ the Lie algebra $\mu_{\gca}\in\sca_{\beta}$ for $\beta\in\tg^+$
conjugate to $\beta_{\mu_{\gca}}$.
\end{proposition}

Recall from graph theory that two edges $l_k,l_m$ of a graph $\gca$ are called {\it
adjacent} if they share a vertex, which will be denoted by $l_k\sim l_m$.  The {\it
line graph} $L(\gca)$ of $\gca$ is the graph whose vertices are the edges of
$\gca$ and where two of them are joined if and only if  they are adjacent.  The {\it
adjacency matrix} $\Adj{\gca}$ of a graph $\gca$ with a labelling $\{ v_1,...,v_p\}$
of the set of vertices is defined as the (symmetric) $p\times p$ matrix with $1$ in
the entry $ij$ if $v_iv_j$ is an edge and zero otherwise.

Since the set $\left\{\alpha_{ij}^{p+k}:l_k=v_iv_j, \; k=1,...,q\right\}$ is
linearly independent, the matrix $U$ associated to $\mu_{\gca}$ is positive definite
and so the linear system
\begin{equation}\label{grafosU}
    U\left[%
\begin{smallmatrix}
  c_1 \\
  \vdots \\
  c_q \\
\end{smallmatrix}%
\right]=\nu\left[%
\begin{smallmatrix}
  1 \\
  \vdots \\
  1 \\
\end{smallmatrix}%
\right]
\end{equation}
admits a unique solution satisfying $\sum\limits_{k=1}^qc_k=1$.  Another way to
define such solution is as a weighting $(c_1,...,c_q)$ on the set of edges $\{
l_1,...,l_q\}$ of $\gca$ such that
\begin{equation}\label{grafosw}
    3c_k+\sum_{l_m\sim l_k}c_m=\nu, \qquad\forall k=1,...,q.
\end{equation}

We note that $\la\alpha_{ij}^{p+k},\alpha_{i'j'}^{p+k'}\ra=1$ if $l_k\sim l_{k'}$
and zero otherwise, from which it follows that $U=3I+\Adj{L(\gca)}$, where
$\Adj{L(\gca)}$ is the adjacency matrix of the line graph of $\gca$.

\begin{definition}\label{defpos}
A graph $\gca$ is said to be {\it positive} if the vector
$$
\left(
3I+\Adj{L(\gca)}\right)^{-1}\left[%
\begin{smallmatrix}
  1 \\
  \vdots \\
  1 \\
\end{smallmatrix}%
\right]
$$
has all its entries positive.
\end{definition}

Thus a graph is positive if and only if the solution $(c_1,...,c_q)$ to
(\ref{grafosU}) (or equivalently the weighting (\ref{grafosw})) satisfies $c_k>0$
for all $k=1,...,q$.

\begin{theorem}\label{grafoseinstein}
The $2$-step nilpotent Lie algebra $\mu_{\gca}$ attached to a graph $\gca$ is an
Einstein nilradical if and only if $\gca$ is positive.
\end{theorem}

\begin{proof}
We may assume that $\gca$ has no isolated vertices since such vertices only
determine an abelian factor of $\mu_{\gca}$ and hence we can apply \cite[Proposition
3.3]{finding}.  Thus the center of $\mu_{\gca}$ coincides with its derived algebra.
If $\mu_{\gca}$ is an Einstein nilradical, that is, there exists
$\mu\in\G.\mu_{\gca}$ such that $S_{\mu}$ is Einstein, then $\Ri$ is conjugate to
$\beta$ and so to $\beta_{\mu_{\gca}}$ (see Proposition \ref{strataEinstein}, (ii)).
It follows from Lemma \ref{2step} that $\Ric_{\mu}$ can never have a zero eigenvalue
(recall that $\mu\simeq\mu_{\gca}$ and so its center and derived algebra coincide),
which implies that $\beta_{\mu_{\gca}}=\sum\limits_{l_k=v_iv_j}c_k\alpha_{ij}^{p+k}$
with $c_1,...,c_q > 0$ since these are precisely the eigenvalues of
$\beta_{\mu_{\gca}}$ with eigenvectors $e_{p+1},...,e_{n}$, respectively. Thus
$(c_1,...,c_q)$ is a solution to (\ref{grafosU}) by Proposition \ref{stratatracy},
(iii), that is, $\gca$ is positive.

Conversely, if $\gca$ is positive then
$\lambda=\mu_{\gca}(\sqrt{c_1},...,\sqrt{c_q})\in\G.\mu_{\gca}$, where
$(c_1,...,c_q)$ is the (positive) solution to (\ref{grafosU}), and $S_{\lambda}$ is
Einstein by Theorem \ref{tracy}.  Thus $\mu_{\gca}$ is an Einstein nilradical.
\end{proof}

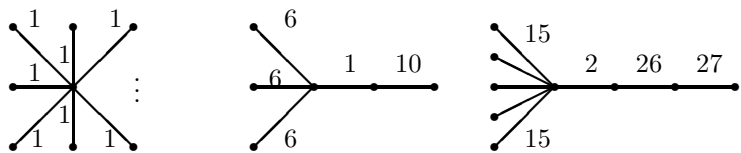
\begin{figure}
\setlength{\unitlength}{.4cm}
\begin{picture}(20,10)
\thicklines

\put(2,4){\circle*{.22}} \put(2,6){\circle*{.22}}\put(2,2){\circle*{.22}}
\put(0,4){\circle*{.22}} \put(4,3.5){\vdots}
\put(0,6){\circle*{.22}}\put(0,2){\circle*{.22}}
\put(4,6){\circle*{.22}}\put(4,2){\circle*{.22}}

\put(2,4){\line(-1,1){2}}\put(2,4){\line(1,1){2}} \put(2,4){\line(-1,-1){2}}
\put(2,4){\line(1,-1){2}} \put(2,4){\line(0,1){2}}
\put(2,4){\line(0,-1){2}}\put(2,4){\line(-1,0){2}}

\put(3.2,6){1}\put(3,2){1} \put(0.5,4.2){1}\put(0.5,6){1}\put(0.6,2){1}
\put(1.5,4.8){1}\put(1.5,2.8){1}

%==================

\put(10,4){\circle*{.22}}\put(14,4){\circle*{.22}}\put(12,4){\circle*{.22}}
\put(8,4){\circle*{.22}}\put(8,6){\circle*{.22}}\put(8,2){\circle*{.22}}

\put(10,4){\line(1,0){4}} \put(8,4){\line(1,0){2}} \put(10,4){\line(-1,1){2}}
\put(10,4){\line(-1,-1){2}}

\put(9,2){6}\put(8.5,4){6}\put(9,6){6} \put(11,4.5){1}\put(12.7,4.5){10}

%==========================

\put(18,4){\circle*{.22}} \put(16,6){\circle*{.22}}\put(16,3){\circle*{.22}}
\put(16,5){\circle*{.22}} \put(16,2){\circle*{.22}}
\put(16,4){\circle*{.22}}\put(20,4){\circle*{.22}}
\put(22,4){\circle*{.22}}\put(24,4){\circle*{.22}}

\put(18,4){\line(1,0){2}} \put(18,4){\line(-1,1){2}}
\put(18,4){\line(-1,-1){2}}\put(18,4){\line(-2,1){2}}
\put(18,4){\line(-2,-1){2}}\put(18,4){\line(-1,0){2}}
\put(20,4){\line(1,0){2}}\put(22,4){\line(1,0){2}}

\put(17,5.5){15} \put(17,2){15} \put(19,4.5){2} \put(20.7,4.5){26}\put(22.7,4.5){27}

\end{picture}\caption{Positive trees with an edge of valency greater than 3.}\label{nn}
\end{figure}

Canonical examples of positive graphs are those for which $L(\gca)$ is regular (i.e.
all the valencies $\val(l_k)=\sharp\{ l_m:l_m\sim l_k\}$ are the same for any $k$).
Indeed, $(1,...,1)$ is a solution to (\ref{grafosw}) if and only if $L(\gca)$ is
regular (note that this is not equivalent to $\gca$ regular, as the graph on the
left in Figure \ref{nn} shows). It is also clear that a graph is positive if and
only if each of its connected components is so. Automorphisms of $L(\gca)$ (i.e.
permutations $\sigma$ of $\{ 1,...,q\}$ such that $l_k\sim l_m$ if and only if
$l_{\sigma(k)}\sim l_{\sigma(m)}$) may be used to facilitate the computation of the
weighting due to the following symmetry property.

\begin{lemma}\label{auto}
Let $(c_1,...,c_q)$ be the weighting defined in {\rm (\ref{grafosw})}.  Then
$c_k=c_{\sigma(k)}$ for any automorphism $\sigma$ of $L(\gca)$ and $k=1,...,q$.
\end{lemma}

\begin{proof}
By using that $l_k\sim l_m$ if and only if $l_{\sigma(k)}\sim l_{\sigma(m)}$ we get
that if $(c_1,...,c_q)$ is a weighting with a given $\nu\in\RR$ then
$(c_{\sigma(1)},...,c_{\sigma(q)})$ satisfies the equations for a weighting with the
same $\nu$ as well, and thus the assertion follows from the uniqueness of the
solution.
\end{proof}

\begin{figure}
\setlength{\unitlength}{.5cm}
\begin{picture}(22,10)
\thicklines

\put(2,3){\circle*{.2}}\put(5,3){\circle*{.2}}
\put(0.5,4.5){\circle*{.2}}\put(0.5,1.5){\circle*{.2}}
\put(6.5,4.5){\circle*{.2}}\put(6.5,1.5){\circle*{.2}} \put(2,3){\line(1,0){3}}
\put(2,3){\line(-1,1){1.5}}\put(5,3){\line(1,1){1.5}} \put(2,3){\line(-1,-1){1.5}}
\put(5,3){\line(1,-1){1.5}} \put(1.5,4){1}\put(1.5,1.5){1}
\put(5.2,4){1}\put(5.2,1.5){1} \put(3.5,2){0}

\put(2.5,0){$\gca_{2,2,0}$}

\put(10,3){\circle*{.2}}\put(14,3){\circle*{.2}}
\put(12,5){\circle*{.2}}\put(12,7){\circle*{.2}} \put(12,9){\circle*{.2}}
\put(10,3){\line(1,0){4}} \put(10,3){\line(1,1){2}}\put(14,3){\line(-1,1){2}}
\put(10,3){\line(1,2){2}}\put(14,3){\line(-1,2){2}}
\put(10,3){\line(1,3){2}}\put(14,3){\line(-1,3){2}} \put(11.1,4.5){1}\put(11.2,6){1}
\put(11,8){1} \put(12.7,4.5){1}\put(12.6,6){1}\put(13,8){1} \put(12,2){0}

\put(10.5,0){$\gca_{0,0,3}$}

\put(18,3){\circle*{.2}}\put(21,3){\circle*{.2}}
\put(16.5,4.5){\circle*{.2}}\put(16.5,1.5){\circle*{.2}} \put(23.5,3){\circle*{.2}}
\put(19.5,6){\circle*{.2}} \put(18,3){\line(1,0){3}}
\put(18,3){\line(-1,1){1.5}}\put(21,3){\line(1,0){2.5}}
\put(18,3){\line(-1,-1){1.5}}
\put(18,3){\line(1,2){1.5}}\put(21,3){\line(-1,2){1.5}} \put(17,4.5){15}
\put(16.5,2.5){15}\put(18.2,4.5){8}\put(20.5,4.5){14}
\put(19.5,2){-1}\put(22,3.5){18}

\put(18.5,0){$\gca_{2,1,1}$}

\end{picture}\caption{Nonpositive graphs and their weightings}\label{patudo}
\end{figure}
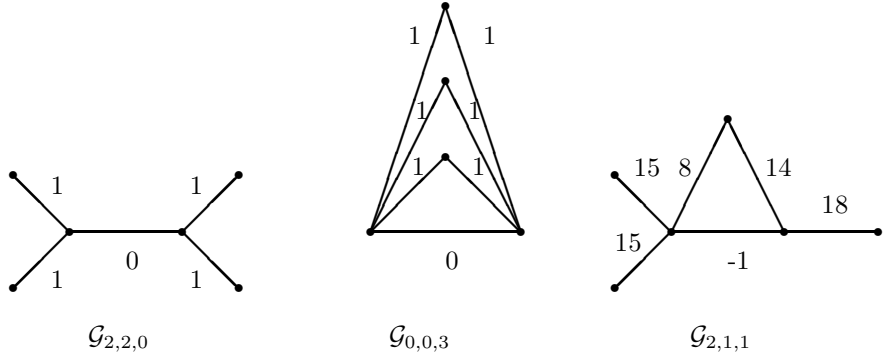

The simplest examples of graphs which are not positive are given in Figure
\ref{patudo}, endowed with their corresponding weightings. Recall that we have
already proved that the $2$-step nilpotent Lie algebra attached to $\gca_{2,2,0}$ is
not an Einstein nilradical in Example \ref{flow2}.  It is not hard to see that any
graph having at most $5$ vertices is positive, with the only exception of
$\gca_{0,0,3}$, the one in the middle in Figure \ref{patudo}.  $\gca_{2,2,0}$ has
$n=p+q=6+5=11$ and gives the lowest dimensional counterexample possible to Problem 4 obtained
by the graph construction.

If $l=vw$ is an edge of a graph then $l$ is said to be {\it incident at} $v$ (or
$w$) and $l$ is called an {\it end edge} if one of the valencies $\val(v), \val(w)$
equals one.

\begin{figure}
\setlength{\unitlength}{.5cm}
\begin{picture}(6,5)
\thicklines \put(2,1){\circle*{.2}}\put(4,1){\circle*{.2}}
\put(0.5,2.5){\circle*{.2}}\put(0.5,-.5){\circle*{.2}}
\put(5.5,2.5){\circle*{.2}}\put(5.5,-.5){\circle*{.2}}
\put(3,2){\circle*{.2}}\put(3,3){\circle*{.2}}

\put(2,1){\line(1,0){2}} \put(2,1){\line(-1,1){1.5}}\put(4,1){\line(1,1){1.5}}
\put(2,1){\line(-1,-1){1.5}} \put(4,1){\line(1,-1){1.5}} \put(0.5,.8){\vdots}
\put(5.3,.8){\vdots}\put(2.9,2.2){\vdots}
\put(2,1){\line(1,1){1}}\put(4,1){\line(-1,1){1}}
\put(2,1){\line(1,2){1}}\put(4,1){\line(-1,2){1}}

\put(0,1){r} \put(6,1){s} \put(3,3.5){t}
\end{picture}\caption{The graph $\gca_{r,s,t}$}\label{patudorst}
\end{figure}
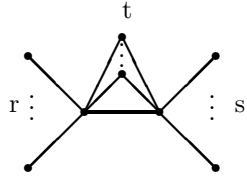

\begin{definition}\label{defgrst}
Let $\gca_{r,s,t}$ be the graph with $p=r+s+t+2$ vertices and $q=r+s+2t+1$ edges
such that there is an edge $l=vw$ which is adjacent to every other edge, $r$ and $s$
are the number of edges incident at $v$ and $w$, respectively, having all other
vertices of valency $1$, and $2t$ is the number of edges adjacent to $l$ but with
all other vertices of valency $2$, as in Figure \ref{patudorst}.  We will always
assume that $r \ge s$ and also that either $s \ne 0$ or $t \ne 0$.  We will say that
a graph $\gca$ {\it contains faithfully} $\gca_{r,s,t}$ if it has an edge $l$ as
above for $r,s,t$ such that the $r+s$ vertices have valency one and the remaining
$t$ vertices have valency two in $\gca$.
\end{definition}

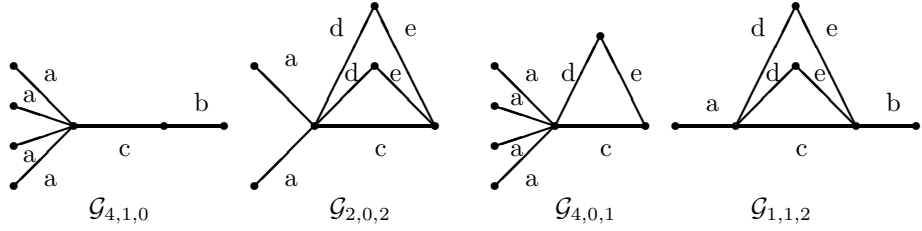
\begin{figure}
\setlength{\unitlength}{.4cm}
\begin{picture}(30,10)
\thicklines

\put(2,4){\circle*{.22}}\put(5,4){\circle*{.22}}\put(7,4){\circle*{.22}}
\put(0,6){\circle*{.22}}\put(0,2){\circle*{.22}}
\put(0,3.35){\circle*{.22}}\put(0,4.65){\circle*{.22}}

\put(2,4){\line(1,0){3}} \put(2,4){\line(-3,1){2}}\put(5,4){\line(1,0){2}}
\put(2,4){\line(-3,-1){2}} \put(2,4){\line(-1,1){2}} \put(2,4){\line(-1,-1){2}}

\put(1,5.5){a}\put(.3,4.8){a} \put(.3,2.8){a}\put(1,2){a}
\put(3.5,3){c}\put(6,4.5){b}

\put(2.5,1){$\gca_{4,1,0}$}

\put(10,4){\circle*{.22}}\put(14,4){\circle*{.22}}\put(12,8){\circle*{.22}}
\put(12,6){\circle*{.22}}\put(8,6){\circle*{.22}}\put(8,2){\circle*{.22}}
\put(10,4){\line(1,0){4}} \put(10,4){\line(1,1){2}}\put(14,4){\line(-1,1){2}}
\put(10,4){\line(1,2){2}}\put(14,4){\line(-1,2){2}} \put(10,4){\line(-1,1){2}}
\put(10,4){\line(-1,-1){2}}
\put(10.5,7){d}\put(11,5.5){d}\put(13,7){e}\put(12.5,5.5){e}
\put(9,6){a}\put(9,2){a} \put(12,3){c}

\put(10.5,1){$\gca_{2,0,2}$}

\put(18,4){\circle*{.22}}\put(21,4){\circle*{.22}}
\put(16,6){\circle*{.22}}\put(16,3.35){\circle*{.22}} \put(16,4.65){\circle*{.22}}
\put(16,2){\circle*{.22}} \put(19.5,7){\circle*{.22}}

\put(18,4){\line(1,0){3}} \put(18,4){\line(-1,1){2}}
\put(18,4){\line(-1,-1){2}}\put(18,4){\line(-3,1){2}} \put(18,4){\line(-3,-1){2}}
\put(18,4){\line(1,2){1.5}}\put(21,4){\line(-1,2){1.5}} \put(17,5.5){a}
\put(17,2){a}\put(16.5,4.7){a} \put(16.5,3){a}\put(18.2,5.5){d}\put(20.5,5.5){e}
\put(19.5,3){c}

\put(18,1){$\gca_{4,0,1}$}

\put(24,4){\circle*{.22}}\put(28,4){\circle*{.22}}\put(26,8){\circle*{.22}}
\put(26,6){\circle*{.22}}\put(22,4){\circle*{.22}}\put(30,4){\circle*{.22}}

\put(24,4){\line(1,0){4}} \put(24,4){\line(1,1){2}}\put(28,4){\line(-1,1){2}}
\put(24,4){\line(1,2){2}}\put(28,4){\line(-1,2){2}} \put(22,4){\line(1,0){2}}
\put(28,4){\line(1,0){2}}

\put(24.5,7){d}\put(25,5.5){d}\put(27,7){e}\put(26.6,5.5){e}
\put(23,4.5){a}\put(29,4.5){b} \put(26,3){c}

\put(24.5,1){$\gca_{1,1,2}$}

\end{picture}\caption{Nonpositive graphs.}\label{patudo2}
\end{figure}

\begin{proposition}\label{np}
Any graph $\gca$ which contains faithfully any of the graphs in Figures
\ref{patudo} and \ref{patudo2} is not positive.
\end{proposition}

\begin{proof}
We first consider a graph $\gca$ which contains faithfully $\gca_{4,1,0}$. Since the
four edges on the left can be interchanged by an automorphism of $\gca$ as they are
end edges, their corresponding weights coincide (see Lemma \ref{auto}).  So if we
call $a,b,c$ to the weights as in Figure \ref{patudo2} then three of the equations
in system (\ref{grafosw}) with $\nu=1$ are
$$
\begin{array}{c}
  6a+c=1-R \\
  3b+c=1-S \\
  4a+b+3c=1-R-S, \\
\end{array}
$$
where $R$ and $S$ denote the sum of the weights of the remaining edges of $\gca$
which are incident at $v$ and $w$, respectively.  By solving this linear system we
get $c=-\tfrac{1}{12}(2R+4S)$.  If we now assume that $\gca$ is positive, then
$R,S\geq 0$ and consequently $c\leq 0$, a contradiction.  We can argue analogously
in the other cases, where we obtain $c=-\unc(R+S+d)$, $c=-\unc(R+S+2d)$ and
$c=-\unc(R+S+\unm(d+e))$ in the remaining cases in Figure \ref{patudo2} from left to
right, and for those in Figure \ref{patudo} we get $c=-\unc(R+S)$, $c=-\unc(R+S)$
and $c=-\unc(R+S+\unm e)$, respectively.  Recall that the triangles can also be
interchanged by an automorphism of $\gca$, and so the weights of the left sides are
all equal to, say, $d$ and those on the right to, say, $e$.
\end{proof}

This gives us the following

\begin{corollary}\label{exceptions}
A graph $\gca_{r,s,t}$ is not positive if and only if one of the following holds:
\begin{itemize}
    \item[(i)] $rs \ge 4$.
    \item[(ii)] $t \ge 3$.
    \item[(iii)] $t \ge 1$ and $t+r/2 \ge 3$.
    \item[(iv)] $r \ge 2$, $s\ge 1$ and $t\ge 1$.
    \item[(v)] $s\ge 1$ and $t\ge 2$.
\end{itemize}
\end{corollary}

\begin{proof}
Any graph $\gca_{r,s,t}$ with $r,s,t$ satisfying any of these conditions is not
positive by Proposition \ref{np}, and the remaining eleven cases are easily seen to
be positive.
\end{proof}

An equivalent way to state Proposition \ref{np} is as follows: any graph faithfully
containing $\gca_{r,s,t}$ such that $r,s,t$ satisfy any of the conditions (i)-(v) in
Corollary \ref{exceptions} is not positive.  This provides a great deal of
counterexamples to Problem 4 in the $2$-step nilpotent case, covering all dimensions
starting from $11$.

In order to obtain an existence result, we are going to study the special case
when $\gca$ is a tree (i.e. a connected graph with no cycles).  We already know that
not every tree is positive (see Figures \ref{patudo} and \ref{patudo2}), but we
shall give a simple sufficient condition.  Let us assume from now on that $\gca$ is
a tree with $q-1$ vertices and $q$ edges such that $\val(l) \le 3$ for any edge $l$
in $L(\gca)$. It is easy to see that for $q \le 5$ such graphs are
\begin{itemize}
\item[(i)] $A_p$, $p=2, \dots, 6$,

\item[(ii)] $D_p$, $p=4, \dots, 6$,

\item[(iii)] $E_6$,

\item[(iv)] $\hca$,
\end{itemize}
where $A_p$, $D_p$ and $E_6$ denote the well known Dynkin diagrams and $\hca$ is the
graph of the family in the left of Figure \ref{nn}, having $4$ edges.  Note that
$\hca$ is special in the sense that, under the hypothesis we have on $\gca$, if
$\hca$ is an induced subgraph of $\gca$ then $\gca = \hca$.  Moreover, $\hca$ is the
only of these graphs having end edges of valency $3$.  It is easy to prove that all
the graphs in (i)-(iv) are positive.

\begin{theorem}\label{exi}
If $\gca$ is a tree with $\val(l) \le 3$ for each edge $l$, then $\gca$ is positive.
\end{theorem}

\begin{proof}
We will prove this result by induction on the number of edges $q$.  By the above
discussion, we may assume from now on that $q>5$ and our inductive hypothesis will
be that any tree with $q-1$ or a smaller number of edges and such that the valencies of the edges
are all at most $3$ is positive.  Given a graph $\gca$ with edges $\{
l_1,...,l_q\}$, for $1\le n \le
 m \le q$, let $U^{(n,m)}$ denote the $(m-n+1) \times (m-n+1)$ minor of $U$ corresponding
to rows and columns from $n$ to $m$.  Also, let $M^{(n,m)}_j$ denote the matrix
obtained by replacing in $U^{(n,m)}$ the $j$-th column by the column vector
$(1,\dots, 1)$.

By Cramer's rule, the $j$-th entries of the solution to the system (\ref{grafosU})
for $U=U^{(n,m)}$ and $\nu=1$ (or equivalently, the $j$-th weight corresponding to
the graph with edges $\{ l_n,...,l_m\}$) is given by
\begin{equation}\label{cramer}
\frac{\det M^{(n,m)}_j}{\det U^{(n,m)}}.
\end{equation}
Recall that $U^{(n,m)}$ is positive definite, thus the sign of the $j$-th weight
equals the sign of $\det M^{(n,m)}_j$.

Let $\gca$ be a tree as in the statement and let us consider an edge $l$ of $\gca$
which is not an end edge. Since $\gca$ is a tree there is a labelling such that if
$l=l_j$ then $l_i \nsim l_k$ for all $i < j < k$. Moreover we can also assume that:
if $\val(l)=2$ then we have $l_{j-1}\sim l_j\sim l_{j+1}$; if $\val(l)=3$ then
$l_{j-1}\sim l_j\sim l_{j+1},l_{j+2}$.  By computing the determinant by the $j$-th
row, a straightforward calculation shows that
\begin{equation}\label{det3}
\begin{array}{rl}
\det M_j^{(1,q)} =& -\det M_{j-1}^{(1,j-1)} \det U^{(j+1,q)}+\det U^{(1,j-1)} \det U^{(j+1,q)}  \\
&- \det M_1^{(j+1,q)} \det U^{(1,j-1)}- u_{j,j+2}\det M_2^{(j+1,q)} \det
U^{(1,j-1)}.
\end{array}
\end{equation}
where $u_{j,j+2}$ equals $1$ or $0$ according to the valency of $l_j$.

We then have two subgraphs of $\gca$ to consider: $\gca_1$ corresponding to the
edges $\{l_1,\dots, l_{j-1}\},$ and $\gca_2$ corresponding to $\{l_{j+1},\dots,
l_q\}$. Let $(a_1,...,a_{j-1})$ and $(b_{j+1},...,b_q)$ denote the weightings of
these graphs corresponding to $\nu=1$. Recall that by inductive hypothesis $a_k>0$
and $b_j>0$ for any $k,j$. Moreover, by (\ref{cramer}) the above equation can be
stated as
\begin{equation}\label{det4}
\begin{array}{r}
\det M_j^{(1,q)} =\det U^{(1,j-1)} \det U^{(j+1,q)} \left[(\frac{1}{2} -
a_{j-1})+(\frac{1}{2} - b_{j+1}-u_{j,j+2}b_{j+2})\right].
\end{array}
\end{equation}

By definition, these weightings satisfy
\begin{equation}
\begin{array}{rr}
(a)& 3a_{j-1}+ u_{j-1,j-2}a_{j-2}+u_{j-1,j-3}a_{j-3} = 1, \\ \\
(b) & 3b_{j+1}+ u_{j+1,j+2}b_{j+2}+u_{j+1,j+3}b_{j+3} = 1,
\end{array}
\end{equation}
where we have included the $u_{k,l}$ factor to consider all possible cases as for
example when $l_{j-1}$ is an end edge. It is easy to see that $(a)$ implies that
$$
\unm-a_{j-1}=\unm(a_{j-1}+u_{j-1,j-2}a_{j-2}+u_{j-1,j-3}a_{j-3})>0,
$$
independently of $\val(l_{j-1})$.

To see that the other term in the right hand side of equation (\ref{det4}) is
nonnegative, suppose first that $\val(l_j)=2 $ (i.e. $u_{j,j+2}=0$). In this case,
by $(b)$ we have that
$$\unm-b_{j+1}=\unm(b_{j+1}+u_{j+1,j+2}b_{j+2}+u_{j+1,j+3}b_{j+3})>0.$$ On the other
hand, if $\val(l_j)=3 $ (i.e. $u_{j,j+2}=1$), we will also consider the equation
corresponding to $l_{j+2}$ (in $\gca_2$),
 $$\begin{array}{rr} (c) & 3b_{j+2}+ b_{j+1}+u_{j+2,j+4}b_{j+4} = 1.
\end{array}$$
Since in this case $u_{j+1,j+2}=1,$ from $(b)$ and $(c)$ we now get
$$\unm-b_{j+1}-b_{j+2}= \unc(u_{j+1,j+3}b_{j+3}+ u_{j+2,j+4}b_{j+4})\ge 0.$$
 In either case from (\ref{det4}) we get that $\det
M_j^{(1,q)}$ is positive, as was to be shown.

Now, if $l$ is an end edge, we can use the same proof considering just $\gca_2$
since we have already excluded the case $\gca=\hca$ and therefore $\val(l)$ is at
most $2$.
\end{proof}

The examples in Figure \ref{nn} show that the condition in Theorem \ref{exi} is not
necessary for a tree to be positive.

Clearly, if the graphs $\gca,\gca'$ are isomorphic then $\mu_{\gca}$ is isomorphic
to $\mu_{\gca'}$ as Lie algebras.  The converse assertion has been proved by M.
Mainkar \cite{Mnk}.  Geometrically, this means that two solvmanifolds
$S_{\mu_{\gca}}$, $S_{\mu_{\gca'}}$ are isometric if and only if the graphs $\gca$,
$\gca'$ are isomorphic.  In particular, Theorem \ref{exi} provides a method to construct a great deal of
examples of nonisometric Einstein solvmanifolds.


\begin{thebibliography}{MM}

\bibitem[B]{Bss} {\sc A. Besse}, Einstein manifolds, {\it Ergeb. Math.} {\bf 10} (1987), Springer-Verlag,
Berlin-Heidelberg.

\bibitem[DM]{DnMnk} {\sc S.G. Dani, M. Mainkar}, Anosov automorphisms on compact nilmanifolds associated with grphs,
{\it Trans. Amer. Math. Soc.} {\bf } (2004).

\bibitem[GH]{GzHkm} {\sc M. Goze, Y. Hakimjanov}, Sur le algebres de Lie
nilpotentes admettant un tore de derivations, {\it Manusc. Math.} {\bf 84} (1994),
115-224.

\bibitem[H]{Hbr} {\sc J. Heber}, Noncompact homogeneous Einstein spaces, {\it Invent. math}. {\bf 133} (1998), 279-352.

\bibitem[K]{Krw1} {\sc F. Kirwan}, Cohomology of quotients in symplectic and algebraic
geometry, {\it Mathematical Notes} {\bf 31} (1984), Princeton Univ. Press,
Princeton.

\bibitem[KMP]{Krd} {\sc K. Kurdyka, T. Mostowski, A. Parusi\'nski}, Proof of the gradient conjecture of R. Thom., {\it
Ann. of Math.} (2) {\bf 152} (2000), 763-792.

\bibitem[L1]{soliton}  {\sc J. Lauret}, Ricci soliton homogeneous nilmanifolds,
{\it Math. Annalen} \textbf{319} (2001), 715-733.

\bibitem[L2]{critical} \bysame, Standard Einstein solvmanifolds as critical points,
{\it Quart. J. Math.} {\bf 52} (2001), 463-470.

\bibitem[L3]{finding}  \bysame, Finding Einstein solvmanifolds by a variational method,
{\it Math. Z.} {\bf 241} (2002), 83-99.

\bibitem[L4]{minimal} \bysame, A canonical compatible metric for geometric structures on
nilmanifolds, {\it Ann. Global Anal. Geom.} {\bf 30} (2006), 107-138.

\bibitem[L5]{praga}  \bysame, Minimal metrics on nilmanifolds, Diff. Geom. and its Appl.,
Proc. Conf. Prague September 2004 (2005), 77-94 ({\it arXiv:} math.DG/0411257).

\bibitem[L6]{standard} \bysame, Einstein solvmanifolds are standard, {\it Ann. of Math.}, in press,
{\it arXiv:} math.DG/0703472.

\bibitem[Mg]{Mgn} {\sc L. Magnin}, Sur les algebres de Lie nilpotents de dimension
$\leq$ 7, {\it J. Geom. Phys.}, 1986 Vol. {\bf III}, 119-144.

\bibitem[Mn]{Mnk} {\sc M. Mainkar}, personal communication, 2006.

\bibitem[M]{Mrn} {\sc A. Marian}, On the real moment map, {\it Math. Res. Lett.} {\bf 8} (2001), 779-788.

\bibitem[Ml]{Mll} {\sc D. Millionschikov}, Graded filiform Lie algebras and symplectic
nilmanifolds, Advances in the Mathematical Sciences (AMS) {\bf 55} (2004), 259-279
({\it arXiv:} math.DG/0205042).

\bibitem[Ms]{Mss} {\sc R. Moussu}, Sur la dynamique des gradients. Existence de vari\'et\'es invariantes, {\it Math. Ann.} {\bf 307} (1997), 445-460.

\bibitem[N]{Nss} {\sc L. Ness}, A stratification of the null cone via the momentum map, {\it Amer. J. Math.} {\bf 106} (1984), 1281-1329 (with an appendix by D. Mumford).

\bibitem[NN]{NktNkn} {\sc E.V. Nikitenko, Yu.G. Nikonorov}, Six-dimensional Einstein solvmanifolds, {\it Siberian
Adv. Math.} {\bf 16} (2006), 66-112.

\bibitem[Nk1]{Nkl1} {\sc Y. Nikolayevsky}, Einstein solvmanifolds with free nilradical, {\it Ann. Global Anal. Geom.} {\bf 33} (2008), 71-87.

\bibitem[Nk2]{Nkl2} \bysame, Nilradicals of Einstein solvmanifolds, {\it arXiv:} math.DG/0612117.

\bibitem[P]{Pyn} {\sc T. Payne}, The existence of soliton metrics for nilpotent Lie
groups, {\it Geom. Ded.} {\bf 145} (2010), 71-88.

\bibitem[S]{Sjm} {\sc R. Sjamaar}, Convexity properties of the moment mapping
re-examined, {\it Adv. Math.} {\bf 138} (1998), 46-91.

\bibitem[W]{Wll} {\sc C.E. Will}, Rank-one Einstein solvmanifolds of dimension $7$,
{\it Diff. Geom. Appl.} {\bf 19} (2003), 307-318.

\end{thebibliography}
\end{document}